\documentclass[11pt]{amsart}
\usepackage{amsfonts}
\usepackage{amsmath}
\usepackage{fancyhdr,amssymb}
\usepackage{amsthm}
\usepackage{ color, ulem}
\input xy
\xyoption {all}

\def \generall {{Z_{X, \, 1,\,\beta\,}(\alpha_1, \ldots, \alpha_{\ell}\,|\,k_1, \ldots, k_{\ell})}}

\newcommand{\comment}[1]{}
\newtheorem{theorem}{Theorem}
\newtheorem {lemma}{Lemma}
\newtheorem{conjecture}{Conjecture}

\newtheorem {proposition}{Proposition}
\theoremstyle{definition}
 
\theoremstyle {definition}

\newcommand{\OO}{\mathcal O}
\newcommand{\BB}{{\mathcal{B}}}
\newcommand{\Ext}{{\text{Ext}}}

\begin{document}
\baselineskip=16pt
\title [Rationality of descendent series for Hilbert and Quot schemes]{Rationality of descendent
  series for Hilbert and Quot schemes of surfaces}
  
\date{March 2021}

\author{Drew Johnson}
\address{Department of Mathematics, ETH Z\"urich}
\email {jared.johnson@math.ethz.ch}
\email {d.johnson@mathematics.byu.edu}
\author{Dragos Oprea}
\address{Department of Mathematics, University of California, San Diego}
\email {doprea@math.ucsd.edu}
\author{Rahul Pandharipande}
\address{Department of Mathematics, ETH Z\"urich}
\email {rahul@math.ethz.ch}

\begin{abstract}
  Quot schemes of quotients of a trivial bundle of arbitrary rank 
  on a nonsingular projective
  surface $X$ carry perfect obstruction theories and virtual fundamental classes
  whenever the quotient sheaf has at most  1-dimensional support.
  The associated generating series of virtual Euler characteristics was conjectured
  to be a rational function in \cite{OP1} when $X$ is simply connected.
  We conjecture here the rationality of
  more general descendent series with insertions obtained from the Chern
  characters of the tautological sheaf.
  We prove the rationality of  descendent series in Hilbert scheme cases
  for all curve classes and
  in  Quot scheme cases when the curve class is $0$.
\end{abstract}

\maketitle

\setcounter{tocdepth}{1}
\tableofcontents

\setcounter{section}{-1}
\section{Introduction}
\subsection{Motivation} 
Let $X$ be a nonsingular
projective surface, and let $X^{[n]}$ denote the Hilbert scheme of points. A well-known formula of G\"ottsche \cite{Go} expresses the topological Euler characteristics of the Hilbert schemes in terms of the Dedekind eta function \begin{equation}\label{gff}\sum_{n=0}^{\infty} \mathsf {e}(X^{[n]})\, q^n=\left(q^{-\frac{1}{24}}\, \eta(q)\right)^{-\mathsf e(X)}.\end{equation} G\"ottsche's formula reflects the action of the Heisenberg algebra on the cohomology of $X^{[n]}$ constructed by \cite{Gr, Na}. 

There are at least two possible directions of extending \eqref{gff}. First, we may view $X^{[n]}$ as the moduli space of rank $1$ sheaves with trivial determinant. The higher rank moduli spaces of sheaves over $X$ play a central role in Vafa-Witten theory \cite{TT, VW}. Explicit expressions for the generating series of the rank $2$ and $3$ moduli spaces are conjectured in \cite{GK1, GK2}. Since the higher rank moduli spaces may be singular, the Euler characteristics are replaced by virtual analogues which take into account the deformation theory of the moduli space. 

In a different direction, we may promote $X^{[n]}$ to more general Hilbert and Quot schemes and
study the corresponding virtual invariants. 

\subsection{Quot schemes and virtual Euler characteristics}
\label{vech}
Let $X$ be a nonsingular
projective surface, let $\beta\in H_2(X,\mathbb{Z})$ be an effective
curve class of $X$, and let $N\geq 1$ be an integer.
Consider the Quot scheme $\mathsf {Quot}_X(\mathbb C^N, \beta, n)$
parameterizing short exact sequences
\begin{equation}\label{qqqq}
  0\to S\to \mathbb C^N\otimes \mathcal O_X\to Q\to 0
  \end{equation}
where $$\text{rank }Q=0\, , \,\,\, c_1(Q)=\beta\, ,\,\,\, \chi(Q)=n\, .$$
As explained in \cite{MOP, OP1},
$\mathsf {Quot}_X(\mathbb C^N, \beta, n)$ carries a canonical
$2$-term perfect obstruction theory and a virtual fundamental class of dimension $$\text{vdim}=\chi(S, Q)=Nn+\beta^2.$$
The virtual fundamental class of the Quot scheme was used in \cite{MOP} to prove Lehn's 
conjecture \cite{Le} for $K3$ surfaces.{\footnote{See \cite{MOP2,MOP3,V} for further
    developments.}}

The virtual Euler characteristic is defined using the virtual tangent complex of the canonical 
obstruction theory \cite{FG}. 
By analogy with the Poincar\'e-Hopf theorem, we set $$\mathsf e^{\text{vir}} (\mathsf {Quot}_X(\mathbb C^N, \beta, n))=\int_{\left[\mathsf {Quot}_X(\mathbb C^N, \beta, n)\right]^{\text{vir}}} c(T^{\text{vir}} \mathsf {Quot}) \ \in \mathbb{Z}\, ,$$
where $c$ denotes the total Chern class. The virtual tangent bundle is given by $$T^{\text{vir}} \mathsf {Quot}=\text{Ext}_X^{\bullet}(S, Q)$$ at each short exact sequence \eqref{qqqq}. 

The generating series of virtual Euler characteristics,
\begin{equation}\label{rrrr}
  Z_{X, N, \beta}=\sum_{n\in \mathbb{Z}} \mathsf e^{\text{vir}} (\mathsf {Quot}_X(\mathbb C^N, \beta, n))\,q^n\,,
\end{equation}
was introduced and studied in \cite {OP1}. For fixed $X$, $N$, and $\beta$,
    $$\text{vdim }\mathsf{Quot}_X(\mathbb C^N, \beta, n)=Nn+\beta^2<0$$ for $n$ sufficiently
    negative, hence $Z_{X, N, \beta}$ has a finite polar part. The following rationality property was conjectured{\footnote{The conjecture can also be made for
        surfaces which are not simply connected, but we will not study
        non simply connected surfaces here (except in the  $\beta = 0 $ case).}}
in \cite{OP1}.
\begin{conjecture} \label{conj}
  Let $X$ be a nonsingular projective simply connected surface, and let $\beta$ be an effective curve class.
  The series $Z_{X, N, \beta}$ is
  the Laurent expansion of a rational function in $q$. 
\end{conjecture} 

\noindent Conjecture \ref{conj} is known to be true in the following five cases:

\vspace{5pt}
\noindent $\bullet$ For all $N\geq 1$, the series $Z_{X,N,\beta}$ is rational if
\begin{itemize}
\item [(i)] $X$ is {\it any} surface and 
  $\beta=0$  \cite {OP1}, 
\item [(ii)] $X$ is a surface of general type\footnote{Property (ii) is proven in \cite{OP1} for simply connected minimal surfaces of general type with $p_g>0$ and a nonsingular canonical divisor. The assumptions other
    than $p_g>0$ were removed in \cite {L}. A similar analysis was
    done in \cite{L} at the level of $\chi_{-y}$-genera. }
  with $p_g>0$ and $\beta$ is {\it any} effective curve class
  \cite{L, OP1}, 
\item [(iii)] $X$ is an elliptic surface\footnote{Property (iii) is proven in \cite{OP2} for simply connected minimal elliptic surfaces. These assumptions were removed in \cite {L} at the level of $\chi_{-y}$-genera.} with $p_g>0$ \cite {L,OP2}.
\end{itemize}
\noindent $\bullet$ For $N=1$, the series $Z_{X,1,\beta}$ is also rational if
\begin{itemize}
\item [(iv)] $X$ is a blow-up and $\beta$ is a multiple of the exceptional divisor \cite {L,OP1},
\item [(v)] $X$ is a $K3$ surface with reduced invariants for primitive curve classes \cite{OP1}.
\end{itemize} 
Our first result here is a resolution of Conjecture \ref{conj} in case $N=1$.

\begin{theorem} \label{t1} Let $X$ be a nonsingular projective simply connected surface, and let $\beta$ be an effective curve class.
  The generating
  series $Z_{X,1,\beta}$ of virtual Euler characteristics is
  the Laurent expansion of a rational function in $q$.
\end{theorem}

In the $N=1$ case, the Quot scheme $\mathsf{Quot}_X(\mathbb C^1, \beta, n)$
is simply a Hilbert scheme of points and curves in $X$.
Theorem
\ref{t1} is therefore about the virtual Euler characteristics of such Hilbert
schemes of surfaces.
A crucial idea in our proof is to transform the geometry to
the moduli space of stable pairs \cite{PT1,PT2} on surfaces and to use the associated
Jacobian fibration.

\subsection{Rationality of descendent series}

How special is the rationality of the generating series \eqref{rrrr} of virtual Euler characteristics?
We propose here a wider rationality statement for descendent series.

Let $X$ be a nonsingular projective simply connected
surface, and let  $\mathsf {Quot}_X(\mathbb C^N, \beta, n)$ 
be the Quot scheme parameterizing quotients \eqref{qqqq}.
Let
\begin{eqnarray*}
\pi_1: \mathsf {Quot}_X(\mathbb C^N, \beta, n) \times X &\rightarrow&  \mathsf {Quot}_X(\mathbb C^N, \beta, n)\, , \\
\pi_2: \mathsf {Quot}_X(\mathbb C^N, \beta, n) \times X &\rightarrow & X\, 
\end{eqnarray*}
be the two projections. 
Let $$\mathcal{Q}\to\mathsf {Quot}_X(\mathbb C^N, \beta, n) \times X$$ be the universal quotient.
For a $K$-theory class $\alpha \in K^0(X)$, we define
$$\alpha^{[n]} = {\mathbf R}\pi_{1*}(\mathcal{Q} \otimes \pi_2^*\alpha) \in K^0(\mathsf {Quot}_X(\mathbb C^N, \beta, n))\, .$$

A generalization of the  series \eqref{rrrr} of virtual Euler characteristics is defined
as follows. Let $\alpha_1,\ldots, \alpha_\ell \in K^0(X)$, and let
$k_1,\ldots, k_\ell$ be non-negative integers.
Set
\begin{multline}\label{rrrrt}
  Z_{X, N, \beta} (\alpha_1,\ldots, \alpha_\ell \,|\,k_1, \ldots, k_{\ell})= \\
  \sum_{n\in \mathbb{Z}} q^{n} \cdot \int_{\left[\mathsf {Quot}_X(\mathbb C^N, \beta, n)\right]^{\text{vir}}}
  \mathsf{ch}_{k_1}(\alpha_1^{[n]}) \cdots \mathsf{ch}_{k_\ell}(\alpha_\ell^{[n]})\, c(T^{\text{vir}}\mathsf {Quot})
  \,.
\end{multline}
The Chern characters in \eqref{rrrrt} may be viewed as 
{\it descendent insertions}. Hence, 
we view $Z_{X, \,N, \,\beta}(\alpha_1, \ldots, \alpha_\ell \,|\,k_1, \ldots, k_{\ell})$
as a  {\it descendent series}.

\begin{conjecture} \label{tttt} The descendent
  series $Z_{X,N,\beta}(\alpha_1,, \ldots, \alpha_\ell \,|\,k_1, \ldots, k_{\ell})$
 is the Laurent expansion of a rational function in $q$.
\end{conjecture}

\noindent We can prove Conjecture \ref{tttt} in case either $\beta=0$
or $N=1$. 

\begin{theorem}\label{t2} Let $X$ be a nonsingular projective surface.
 For $\beta=0$, the series
  $$Z_{X,\,N,\,0}(\alpha_1,\ldots, \alpha_\ell \,|\,k_1, \ldots, k_{\ell})\in \mathbb{Q}((q))$$
  is the Laurent expansion of a rational function in $q$.
\end{theorem}

\begin{theorem} \label{t3} Let $X$ be a nonsingular projective simply connected surface, and let $\beta$ be an effective curve class.
  For $N=1$, the 
  series
  $$Z_{X,\,1,\,\beta}(\alpha_1, \ldots, \alpha_\ell \,|\,k_1, \ldots, k_{\ell})\in \mathbb{Q}((q))$$
is
  the Laurent expansion of a rational function in $q$.
\end{theorem}

The rationality statements for surfaces here are parallel to the rationality
of the descendent series for stable pairs on $3$-folds, see \cite{P}
for a survey and 
\cite{PP1,PP2,PP3,PT1,PT2} for foundational results.
Whether the descendent series \eqref{rrrrt} satisfy relations such as
the Virasoro constraints for stable pairs \cite{OOP, MOOP} is an interesting question for further study.

Descendent integrals against the
(non-virtual) fundamental class of the Hilbert scheme of points of a surface
have been studied by Carlsson  \cite{C}; the descendent series are proven
to be quasi-modular. The virtual fundamental class regularizes the
descendent geometry in two ways: the theory can be defined more generally for
Quot schemes of quotients supported on curves  {\it and} the answers are rational functions.

The study of the virtual invariants of Quot schemes of surfaces can also be considered in $K$-theory. For recent results and conjectures related to the rationality of descendent series in $K$-theory, see \cite{AJLOP}. Section $3.2$ of \cite{AJLOP} also contains subsequent developments regarding Conjecture \ref{tttt} for all surfaces with $p_g>0$. 

\subsection {Acknowledgments} 
Our study of the virtual Euler characteristics of the
Quot scheme of surfaces was motivated in part by the
Euler characteristic 
calculations of L. G\"ottsche and M. Kool \cite{GK1, GK2} for the
moduli spaces of rank $2$ and $3$ stable sheaves on surfaces.
We thank A. Marian, W. Lim, A. Oblomkov, A. Okounkov, and R. Thomas for
related discussions.

D. J. was supported by SNF-200020-182181.
D. O. was supported by the NSF through grant DMS 1802228. R.P. was supported by the Swiss National Science Foundation and
the European Research Council through
grants SNF-200020-182181, 
ERC-2017-AdG-786580-MACI, SwissMAP, and the Einstein Stiftung.  
We thank the Shanghai Center for Mathematical Science at Fudan
University for a very productive visit in September 2018 at the start of the
project.

The project has received funding from the European Research
Council (ERC) under the European Union Horizon 2020 Research and
Innovation Program (grant No. 786580).

\section{Virtual Euler characteristics: Theorem \ref{t1}}

\label{vic}

\subsection{Obstruction theory}
We start the proof of Theorem \ref{t1} with an explicit
description of the Hilbert scheme and the obstruction theory
in the $N=1$ case.

Let $X$ be a nonsingular projective surface.
When $N=1$, the following isomorphism was proved in \cite{F}:
$$\mathsf {Quot}_X(\mathbb C^1, \beta, n)\simeq X^{[m]}\times
\mathsf{Hilb}_{\beta}\, .$$ 
Here, $X^{[m]}$ is the Hilbert scheme of $m$ points of $X$,
$\mathsf {Hilb}_{\beta}$ is the Hilbert scheme of divisors of $X$
in the class $\beta$, and $$m=n+\frac{\beta(\beta+K_X)}{2}\, .$$ Under this isomorphism, each pair $(Z, D)\in X^{[m]}\times
\mathsf{Hilb}_{\beta}$ yields a short exact sequence $$0\to I_Z(-D)\to \mathcal O_X\to Q\to 0\, .$$
The Hilbert scheme $\mathsf {Hilb}_{\beta}$ parameterizes
only pure dimension $1$ subschemes. 
There is an Abel-Jacobi map $$\mathsf {AJ}_{\beta}: \mathsf {Hilb}_{\beta}\to \text{Pic}^{\beta} (X)\, ,\,\ \  D\mapsto \mathcal O_X(D)\, ,$$ with fibers
given by projective spaces of possibly varying dimension. As noted in \cite{DKO}, $\mathsf {Hilb}_{\beta}$ carries a virtual fundamental class of dimension $$\text{vdim}_{\beta}=\frac{\beta(\beta-K_X)}{2}\, .$$ 

The virtual fundamental class of $\mathsf {Quot}_X(\mathbb C^1, \beta, n)$ was identified in \cite{L} to equal \begin{equation}\label{wl}\left[\mathsf {Quot}_X(\mathbb C^1, \beta, n)\right]^{\text{vir}}=\mathsf e(\mathsf B)\cap \left(\left[X^{[m]}\right] \times \left[\mathsf{Hilb}_{\beta}\right]^{\text{vir}}\right) \end{equation} where $$\mathsf B=RHom_{\pi}(\mathcal O_{\mathcal W}, \mathcal O(\mathcal D))\, .$$ Here $$\mathcal W\subset X\times  X^{[m]},\,\,\, \mathcal D\subset X\times \mathsf {Hilb}_{\beta}$$ are the universal families, and $$\pi: X\times X^{[m]}\times \mathsf {Hilb}_{\beta}\to  X^{[m]}\times \mathsf {Hilb}_{\beta}$$ is the projection. 

When $X$ is simply connected, the Hilbert scheme $\mathsf {Hilb}_{\beta}=\mathbb P$ is a projective space of dimension $h^0(\beta)-1$. The obstruction bundle for $\mathsf {Quot}_X(\mathbb C^1, \beta, n)$ given above simplifies to the expression found in \cite {OP1}: \begin{equation}\label{obun}\mathsf{Obs} = (H^1(M)-H^0(M))\otimes \mathcal L+\left(M^{[m]}\right)^{\vee}\otimes \mathcal L+ \mathbb C^{p_g}.\end{equation} Here $$M=K_X-\beta$$ and the superscript $(\,\,)^{[m]}$ denotes the usual tautological bundle over the Hilbert scheme of points $X^{[m]}$. Furthermore, $$\mathcal L=\mathcal O_{\mathbb P}(1)\, .$$

Theorem \ref{t1} is established whenever $p_g>0$. For surfaces of positive Kodaira dimension, the claim follows by cases (ii) and (iii)
discussed after Conjecture \ref{conj} in Section \ref{vech}.
The only remaining cases are $K3$ surfaces and their successive
blowups.
Invariants of $K3$ surfaces vanish unless $\beta=n=0$, see \cite {MOP}.
Theorem $6$ of \cite {L} determines the invariants of blowups
in terms of explicit rational functions, see also Section \ref{nms} below.

We assume $p_g=0$ for the remainder of Section \ref{vic}.
Since $\beta$ is an effective curve class, the condition $p_g=0$ implies
$$H^0(M)=H^0(K_X-\beta) = 0\, .$$
The obstruction bundle therefore
further simplifies to $$\mathsf{Obs} = H^1(M)\otimes \mathcal L+\left(M^{[m]}\right)^{\vee}\otimes \mathcal L\, .$$

\subsection{Rationality} \label{rational}For a nonsingular scheme
$S$ endowed with a perfect obstruction theory and obstruction bundle $\mathsf {Obs}$, the virtual Euler characteristic is given by 

$$\mathsf e^{\text{vir}}(S)=\int_{M} \mathsf e(\mathsf {Obs}) \cdot \frac{c(TS)}{c(\mathsf {Obs})}\, .$$
In our situation (assuming $p_g=0$), 
$$\mathsf e^{\text{vir}}(\mathsf {Quot}_X(\mathbb C^1, \beta, n))=\int_{X^{[m]}\times \mathbb P} c_1(\mathcal L)^{h^1(\beta)}\cdot \mathsf e\left(\mathcal L\otimes \left(M^{[m]}\right)^{\vee}\right) \cdot \frac{c(TX^{[m]}) \cdot c(\mathcal L)^{\chi(\beta)}}{c\left(\mathcal L\otimes \left(M^{[m]}\right)^{\vee}\right)}\, .$$ We can integrate out the hyperplane class to reduce the dimension of the projective space to $\chi(\beta)-1$. Theorem \ref{t1} follows from the following result.

\begin{proposition}\label{p1} Let $V$ be a finite dimensional vector space, and let $M\to X$ be a line bundle over a nonsingular
  projective surface. The series $$Z_{X, M, V}=\sum_{n=0}^{\infty} q^n \cdot \int_{X^{[n]}\times \mathbb P(V)} \mathsf e\left(\mathcal L\otimes \left(M^{[n]}\right)^{\vee}\right) \cdot \frac{c(TX^{[n]}) \cdot c(T\mathbb P(V))}{c\left(\mathcal L\otimes \left(M^{[n]}\right)^{\vee}\right)}$$ is a rational function in $q$. 
\end{proposition}

In fact, we will prove a stronger claim.
For a rank $r$ vector bundle $E\to S$ over a scheme $S$ with Chern roots $x_1, \ldots, x_r$, define \begin{equation}\label{pbe}P_d(E)=\sum_{i=1}^{r} \frac{1}{(1+x_i)^d}\, .\end{equation} For a finite sequence $B=(b_1, \ldots, b_\ell)$ of non-negative integers, we set $$P(E, B)=\prod_{i=1}^{\ell} P_{i}(E)^{b_i}.$$ 
Write $$Z_{X, M}[a,B] =\sum_{n=0}^{\infty} q^n\cdot \int_{X^{[n]}}  c_{n-a}\left(\left(M^{[n]}\right)^{\vee}\right) \cdot c(TX^{[n]}) \cdot \frac{P\left(\left(M^{[n]}\right)^{\vee},B\right)}{c\left(\left(M^{[n]}\right)^{\vee}\right)}\, .$$
\begin{proposition}
  \label{p1a}For all pairs $(X, M)$, non-negative integers $a$, and finite sequences $B$, the series $Z_{X,M}[a,B]$ is a rational function in q.
\end{proposition}

Proposition \ref{p1a} implies Proposition \ref{p1} by the
following argument. Let $\zeta=c_1(\mathcal L)$ denote the hyperplane class on $\mathbb P(V)$. We analyze the expressions appearing in Proposition \ref{p1}. First, 
$$\mathsf e\left(\mathcal L\otimes \left(M^{[n]}\right)^{\vee}\right)=\sum_{a=0}^{n} \zeta^a \cdot c_{n-a} \left(\left(M^{[n]}\right)^{\vee}\right)\, .$$ 
Next, we write $x_1, \ldots x_n$ for the Chern roots of $M^{[n]}.$ We have 
$$\frac{1}{c\left(\mathcal L\otimes \left(M^{[n]}\right)^{\vee}\right)} =  \prod_{i=1}^{n} \frac{1}{1 -x_i +\zeta}\, .$$ We expand
$$\frac{1}{1-x_i+\zeta}=\frac{1}{1-x_i}\cdot \sum_{j=0}^{\infty} (-1)^{j}\cdot \zeta^{j} (1 - x_i)^{ -j }$$
which yields
$$\frac{1}{c\left(\mathcal L\otimes \left(M^{[n]}\right)^{\vee}\right)} = \frac{1}{c\left( \left(M^{[n]}\right)^{\vee}\right)} \cdot \left(\sum_{j=0}^{\infty} (-1)^{j} \zeta^{j} H_j \right)\, ,$$ 
where 
$$H_j= \sum_{j_1+\ldots+j_n=j} (1 - x_1)^{ -j_1} \cdots (1 - x_n)^{- j_n}\, .$$
The integral in Proposition \ref{p1} becomes $$\int_{X^{[n]}\times
  \mathbb P(V)}\left(\sum_{a=0}^{n} \zeta^a \cdot c_{n-a} \left(\left(M^{[n]}\right)^{\vee}\right)\right)\cdot \frac{c(TX^{[n]})}{c\left( \left(M^{[n]}\right)^{\vee}\right)} \cdot \left(\sum_{j=0}^{\infty} (-1)^{j} \zeta^{j} H_j \right)\cdot (1+\zeta)^{v}$$ where $\dim V=v$.

After integrating out $\zeta$ over $\mathbb P(V)$,
we are led to expressions of the form 
$$\int_{X^{[n]}}  c_{n-a}\left(\left(M^{[n]}\right)^{\vee}\right) \cdot \frac{c(TX^{[n]})}{c\left( \left(M^{[n]}\right)^{\vee}\right)}\cdot H_j$$ 
with $a+j\leq v-1$. Crucially, both $a$ and $j$ are bounded by $\dim V=v$, independently of $n$. Furthermore, each $H_j$ is symmetric in the Chern roots so can be expressed as a polynomial in the power sums
$$P_d=\sum_{i=1}^{n} \frac{1}{(1-x_i)^{d}}$$
in a fashion which is independent of $n$.
Explicitly, we have $$\sum_{j=0}^{\infty} t^j H_j = \exp \left(\sum_{d=1}^{\infty} \frac{t^d P_d}{d}\right)\, .$$ These remarks reduce the proof of Proposition \ref{p1} to Proposition \ref{p1a}.

\subsection{Proof of Proposition \ref{p1a}}
\subsubsection{Strategy}
 We will prove Proposition \ref{p1a} in two steps:

 \begin{itemize}
 \item [(i)] We first reduce to special rational geometries
   via universality considerations.
\item [(ii)] A geometric argument using the
  moduli space of stable pairs
  will be given for rational surfaces $X$ with a sufficiently positive line
  bundle $M$. 
\end{itemize} 
\subsubsection{Universality}\label{uni}
Fix $\ell\geq 0$. We form the generating series $$Y_{X, M}^{(\ell)}=\sum_{B=(b_1, \ldots, b_\ell)} \frac{z_1^{b_1}}{b_1!} \cdots \frac{z_{\ell}^{b_{\ell}}}{b_{\ell}!} \sum_{n\geq 0} \sum_{a\geq 0}q^n t^a \cdot \int_{X^{[n]}} c_{n-a}\left(\left(M^{[n]}\right)^{\vee}\right) \, c(TX^{[n]}) \, \frac{P\left(\left(M^{[n]}\right)^{\vee},B\right)}{c\left(\left(M^{[n]}\right)^{\vee}\right)}\, .$$ 

\noindent The above expression is multiplicative in the sense that if
$X=X_1\sqcup X_2$, then
\begin{equation}\label{vv44}
  Y_{X, M}^{(\ell)}=Y_{X_1, M_1}^{(\ell)}\cdot Y_{X_2, M_2}^{(\ell)}\, ,
  \end{equation}
  where $M_1, M_2$ are the restrictions of $M$ to $X_1, X_2$ respectively.
Claim \eqref{vv44} is a consequence of the following observations $$X^{[n]}=\bigsqcup_{n_1+n_2=n} X_1^{[n_1]}\times X_2^{[n_2]}$$ $$M^{[n]}=\bigsqcup_{n_1+n_2=n} M_1^{[n_1]}\boxplus M_2^{[n_2]}$$ $$P_{i}\left(\left(M^{[n]}\right)^{\vee}\right)=\bigsqcup_{n_1+n_2=n} P_{i}\left(\left(M_1^{[n_1]}\right)^{\vee}\right)+P_{i}\left(\left(M_2^{[n_2]}\right)^{\vee}\right).$$ The factorials in the definition of $Y_{X, M}^{(\ell)}$ are engineered to offset the prefactors appearing in the binomial expansion $P_i^{b_i}$ of the third identity above. 

As a consequence of above multiplicativity and
the arguments of \cite {EGL}, we have $$Y_{X, M}^{(\ell)}=A_1^{K_X^2}  \cdot A_2^{\chi(\mathcal O_X)} \cdot A_3^{M\cdot K_X} \cdot A_4^{M^2},$$ for universal series $A_1, A_2, A_3, A_4$ in the variables $q, t, z_1,\ldots,z_\ell$.
To prove Proposition \ref{p1a},
we must show that $$\text{Coefficient of }\  t^{a} z_1^{b_1}\cdots z_{\ell}^{b_{\ell}}
\ \  \text{in}\ \  A_1^{K_X^2}  \cdot A_2^{\chi(\mathcal O_X)} \cdot A_3^{M\cdot K_X} \cdot A_4^{M^2}$$ is a rational function in $q$. 

Our method is to study special geometries $(X, M)$. Several choices are possible here\footnote{The simplest geometry $X=\mathbb P^2$ places numerical restrictions leading, at least a priori, to less precise results regarding the denominators of the answers.}, for instance we could pick 
\begin{itemize}
\item [(a)] $X$ is the blowup of $\mathbb P^2$ at $1$ point and $M=dH-eE$, 
\item [(b)] $X$ is the blowup of $\mathbb P^2$ at $2$ points and $M=dH-e_1 E_1 - e_2 E_2$. 
\end{itemize}

For the arguments of the following subsection, we will require
$M$ sufficiently positive. For a concrete discussion,
the results of \cite {R} are useful.
Specifically, if $\kappa$ is a fixed integer, a line bundle $M$, assumed not to equal a multiple of $(-K_X)$, is $\kappa$-very ample provided that the following inequalities hold
\begin{itemize}
\item [(a$'$)] $d\geq e+\kappa, \,\,\,\,e\geq \kappa,$
\item [(b$'$)] $d\geq e_1+e_2+\kappa, \,\,\,\,e_1\geq \kappa, \,\,\,\,e_2\geq \kappa$.
\end {itemize} 
We will furthermore assume\footnote{In the absence of (c), we have less control on the denominators of the rational functions thus obtained.} 
\begin{itemize}
\item [(c)] there exists a divisor $L$ on $X$ such that $L\cdot M=1$.  
\end{itemize} 
Such an $L$ can be chosen in the form $$L=d'H-e'E\,\, \text { or }\,\, L=d'H-e_1'E_1-e_2'E_2'$$ provided 
\begin{itemize} 
\item [(c$'$)] $\gcd (d, e) =1 \text{ and } \gcd (d, e_1, e_2)=1$.
\end{itemize}
To complete the proof of Proposition \ref{p1a} for arbitrary geometries, we
need the following result.

\begin{lemma}
\label{ln}
  Fix $\ell\geq 0$ and $\kappa>0$. Assume that for all $0\leq a\leq \kappa$, and all nonnegative $b_1, \ldots, b_{\ell}$, $${\mathrm{Coefficient\ of}}\ \
  t^{a} z_1^{b_1}\cdots z_{\ell}^{b_{\ell}} \
\ {\mathrm{in}}\ \ A_1^{K_X^2}  \cdot A_2^{\chi(\mathcal O_X)} \cdot A_3^{M\cdot K_X} \cdot A_4^{M^2}$$ is a rational function in $q$ for $(X, M)$ as above. Then the same coefficients are rational in $q$ for all pairs $(X, M).$
\end{lemma}

\proof Examples (a) and (b) give the rationality of the relevant coefficients in the expressions $$A_1^{8}  \cdot A_2 \cdot A_3^{-3d+e} \cdot A_4^{d^2-e^2}\ \text{ and }\ \ A_1^{7}  \cdot A_2 \cdot A_3^{-3d+e_1+e_2} \cdot A_4^{d^2-e_1^2-e_2^2}\, .$$
By varying $d, e, e_1, e_2$ for sufficiently large values with respect to $\kappa$ subject to the conditions above,  we can reconstruct $A_1, A_2, A_3, A_4$ and conclude that their corresponding coefficients are rational in $q$.  
\qed

\subsubsection{Special geometries}\label{uni2} We verify here the hypotheses of Lemma \ref{ln} for pairs $(X, M)$ satisfying all conditions above. The argument however applies more generally for sufficiently positive line bundles $M\to X$. \vskip.1in

To keep the notation simple, we assume  $B=\emptyset$ throughout
Section \ref{uni2}.
Thus \begin{equation}\label{zxma}Z_{X, M}[a]=\sum_{n\geq 0}^{\infty} q^n \cdot \int_{X^{[n]}}  c_{n-a}\left(\left(M^{[n]}\right)^{\vee}\right) \cdot c(TX^{[n]}) \cdot s\left(\left(M^{[n]}\right)^{\vee}\right)\, ,\end{equation}
where $s$ denotes the Segre class. We will indicate how to proceed with the general case $B\neq \emptyset$ in Section \ref{beee}.

We begin by representing the Chern class $c_{n-a}\left(M^{[n]}\right)$ by a natural geometric cycle. To this end, we pick a general linear system $|V|$ in $|M|$ satisfying the following
two properties:
\begin{itemize}
\item [(i)] $\dim |V|=a$, 
\item [(ii)] the curves in $|V|$ are irreducible and reduced. 
\end{itemize}
This can be achieved if the coefficient $d$ of the hyperplane class in $M$ is chosen sufficiently large. Specifically, by \cite[Proposition 5.1]{KT}, the assumption (ii) is satisfied as soon as $M$ is $(2a+1)$-very ample. We write $$\pi:\mathcal C\to |V|$$ for the universal curve. When regarded as the base of $\pi$, we write $\BB$ instead of $|V|$. Let
$$\pi: (\mathcal C/\BB)^{[n]}\to \BB$$ denote the relative Hilbert scheme of points. For all $n$, the space $(\mathcal C/\BB)^{[n]}$ is a nonsingular
projective variety of dimension $$\dim (\mathcal C/\BB)^{[n]}=n+a\, $$ 
by \cite[Theorem 46]{GS}. The assertion uses the assumption that $M$ is sufficiently positive, in particular, we need $M$ to be $a$-very
ample. Furthermore, we have a natural morphism $$j: (\mathcal C/\BB)^{[n]}\to X^{[n]}\, .$$

 Pick $s_0, \ldots, s_{a}$ a basis for $|V|$, viewed as sections of $M$. Each section $s$ of $M$ induces a tautological section $s^{[n]}$ of the bundle $M^{[n]}$ via restriction $$\xi\to s_{\xi}\,
 ,\ \ \,\, s_{\xi}\in H^0(M\otimes \mathcal O_{\xi})=M^{[n]}|_{\xi}\, .$$ Here $\xi\subset X$ is a length $n$ subscheme of $X$. 
 We therefore obtain sections $$s_0^{[n]}, \ldots, s_a^{[n]}$$ of $M^{[n]}\to X^{[n]}.$ The degeneracy locus of these sections consists of subschemes $\xi$ of $X$ such that $$\xi\subset \mathcal C_b$$ for some curve $\mathcal C_b$ of the linear system $|V|$. We therefore conclude
\begin{equation}\label{eqeq}
  j_{\star} (\mathcal C/\BB)^{[n]}=c_{n-a}(M^{[n]})\cap \left[X^{[n]}\right]\, .
 \end{equation}

 We can rewrite \eqref{zxma} using 
equality \eqref{eqeq} as 
\begin{eqnarray*}Z_{X, M}[a]&=&\sum_{n= 0}^{\infty} q^n \cdot \int_{X^{[n]}}  c_{n-a}\left(\left(M^{[n]}\right)^{\vee}\right) \cdot c(TX^{[n]}) \cdot s\left(\left(M^{[n]}\right)^{\vee}\right)\\&=&\sum_{n=0}^{\infty} q^n(-1)^{n-a}                                         \int_{(\mathcal C/\BB)^{[n]}} j^{\star} c(TX^{[n]}) \cdot j^{\star} s\left(\left(M^{[n]}\right)^{\vee}\right)\\&=& (-1)^{a} \,Z_{\mathcal C/\BB,\, M} 
(-q)\, ,\end{eqnarray*} where we define $$Z_{\mathcal C/\BB,\, M}(q)=\sum_{n=0}^{\infty}q^n \int_{(\mathcal C/\BB)^{[n]}} j^{\star} c\left(TX^{[n]}-\left(M^{[n]}\right)^{\vee}\right)\, .$$ 
\vskip.1in                     
We prove the rationality of $Z_{\mathcal C/\BB, \,M}$. The key step is to show that the generating series encodes expressions of the form 
\begin{equation}\label{special}
(-1)^n \left(p_1(n)+2^n \cdot p_2(n)\right)\,\,\,\text {for polynomials }p_1, p_2\, .
\tag{$\star$}
\end{equation}
Series of the form $$\sum_{n=0}^{\infty} (-1)^n \left(p_1(n) +2^n\cdot p_2(n)\right) q^n$$ are rational functions in $q$.\footnote{As a consequence, the denominators of the series of Euler characteristics \eqref{rrrr} are products of $1-q$ and $1-2q$ with various exponents. The same assertion holds true for the
  descendent series of Theorem \ref{t3}. The example of Subsection \ref{r8} with  $\beta=0$ also has
  the same denominators.}
Hence, we will deduce Proposition \ref{p1a} from the following result.

\begin{lemma}\label{l1}
  For sufficiently positive line bundles $M\to X$ satisfying conditions (a$'$), (b$'$), and (c$'$), and families of curves $\mathcal C\to \mathcal B$ satisfying (i) and (ii), the expression 
  \begin{equation}\label{e1} \int_{(\mathcal C/\BB)^{[n]}} j^{\star} c\left(TX^{[n]}-\left(M^{[n]}\right)^{\vee}\right)\end{equation} is of the form \eqref{special}
  for polynomials $p_1(n)$ and $p_2(n)$.
\end{lemma}

\subsubsection{Proof of Lemma \ref{l1}} \label{L1L1L1}
We let $\mathcal H\to \mathcal C$ denote a relatively ample line bundle for the family $$\pi:\mathcal C\to \BB.$$ For instance, we may pick $$\mathcal H=j^{\star} L$$ for the line bundle $L$ whose existence was assumed in (c).
Then, $\mathcal H$ has fiber degree $1$. \vskip.1in

The following structures will play an important role in the proof of
Lemma \ref{l1}: 
\begin{itemize}
\item [(i)] the relative moduli space $\mathfrak M\to \mathcal B$ of torsion free rank $1$ sheaves of degree $0$ over the fibers of $\pi:\mathcal C\to \BB$,
\item [(ii)] the universal sheaf $$\mathcal J\to \mathfrak M\times_{\BB}\mathcal C\, $$ constructed in \cite {AK} for families of reduced irreducible curves, 
\item [(iii)] the universal subscheme $$\mathcal Z_n\hookrightarrow (\mathcal C/\BB)^{[n]}\times_{\BB}\mathcal C\,,$$ 
\item [(iv)] the universal subscheme $\mathcal W_n$ of $X^{[n]}\times X$. 
\end{itemize}    

We write $$\widehat \pi: \mathfrak M\times_{\BB}\mathcal C\to \mathfrak M$$ for the base change of $\pi:\mathcal C\to \BB$.
We consider the sheaves $$\mathcal J\, , \mathcal H\to
\mathfrak M\times_{\BB}\mathcal C$$
where pullback from $\mathcal C$ is understood for the second line bundle. We set $$p_n: \mathbb P_n=\mathbb P\left(\widehat \pi_{\star} \left(\mathcal J\otimes \mathcal H^{n}\right)\right)\to \mathfrak M\, .$$ For $n$ sufficiently large, $\mathbb P_n$ has fibers of constant dimension (by cohomology
vanishing), so $\mathbb P_n$ a projective bundle over $\mathfrak M$.
We write $$\zeta_n=\mathcal O_{\mathbb P_n}(1)\, .$$ 

We will regard the
relative Hilbert scheme $(\mathcal C/\mathcal B)^{[n]}$ as a (subspace of the) moduli space $\mathcal P$ of stable pairs
$$(F, s: \mathcal O_X\to F)$$
on $X$ as explained in \cite[Proposition B8]{PT2}. Here, $$c_1(F)=c_1(M), \quad \chi(F)=1-g+n,$$ with $g$ denoting the arithmetic genus of the linear series $|M|$. We furthermore require that the support of $F$ be contained in $\mathcal B=|V|$. The correspondence between the relative Hilbert scheme and stable pairs can be summarized as follows. For each subscheme $$\xi\subset \mathcal C_b\, ,$$ the canonical sequence $$0\to I_\xi\to \mathcal O_{\mathcal C_b}\to \mathcal O_{\xi}\to 0$$ dualizes to $$0\to \mathcal O_{\mathcal C_b}\to Hom_{\mathcal C_b}(I_{\xi}, \mathcal O)\to \mathcal {E}xt^1_{\mathcal C_b}(\mathcal O_{\xi}, \mathcal O)\to 0\, ,$$ where the last term has dimension zero and length $n$.
Setting $$F= I_\xi^{\vee}=Hom_{\mathcal C_b}(I_{\xi}, \mathcal O),$$
we obtain a stable pair $$s:\mathcal O_X\to F$$ on $X$ with the stated numerical invariants. By a result of \cite{PT2},
$$Ext^{\geq 1}_{\mathcal C_b}(I_{\xi}, \mathcal O)=0\, .$$ Hence,
the above dual can be interpreted as $RHom_{\mathcal C_b}(I_{\xi}, \mathcal O)$ in the derived category. 

As a consequence of the above identifications, there is a natural morphism \begin{equation}\label{taun}
  \tau_n:(\mathcal C/\BB)^{[n]}\to \mathbb P_n\, ,\end{equation}
Indeed, for the moduli space of stable pairs,
we have a natural morphism
\begin{equation}
  \label{ff55ff}
  \mathcal P\to \mathfrak M\, ,\ \ \, \, (F, s:\mathcal O_X\to F)\mapsto F\otimes \mathcal H^{-n}\, .
  \end{equation}
We used here that $\mathcal H$ has fiber degree $1$, so that the twist $F\otimes \mathcal H^{-n}$ has fiber degree $0$.
The fiber of the morphism \eqref{ff55ff}
over a sheaf $J\in \mathfrak M$ is $$\mathbb PH^0(J\otimes \mathcal H^{n})\, .$$

The universal structure $$\mathcal Z_n \, \hookrightarrow\,  (\mathcal C/\BB)^{[n]}\times_{\BB}\mathcal C\, \to\,  \mathbb P_n\times_{\BB}\mathcal C\, \to\,  \mathfrak M\times_{\BB} \mathcal C$$ satisfies
    \begin{equation}\label{ideal}\mathcal I_{\mathcal Z_n}^{\vee}=\mathcal J\otimes
  \mathcal H^{n}\otimes \zeta_n\, .\end{equation} In the above, duals are interpreted in the derived category. 

We now examine the integrand which appears in Lemma \ref{l1}. The following tautological structures over $\mathfrak M$ will be needed in the analysis. 
\begin{itemize}
\item [(A)] Consider the diagram \begin{center}
$\xymatrix{& \mathcal C  \ar[d]^{\pi}  \ar[r]^{j} & X\,. \\ \mathfrak M\ar[r]^{p} & \BB}$
\end{center}
For a bundle $W\to X$, we define
$$\overline W=p^{\star} {\mathbf R}\pi_{\star} j^{\star} W\to \mathfrak M\, .$$ 
\item [(B)] Consider the diagram \begin{center}
$\xymatrix{\mathcal C\times_{\BB} \mathfrak M  \ar[d]^{\widehat \pi}  \ar[r] & \mathcal C\, . \\ \mathfrak M}$ 
\end{center} For a bundle $\mathcal V\to \mathcal C$, we set
\begin{eqnarray*}
  \mathcal V_n\to \mathfrak M, \,\,\,\, & & \mathcal V_n={\bf R}\widehat \pi_{\star}(\mathcal V\otimes \mathcal J^{\vee}\otimes \mathcal H^{-n})\, ,\\
\mathcal V'_n\to \mathfrak M, \,\,\,\, & & \mathcal V'_n=\text{Ext}^{\bullet}_{\,\widehat \pi} (\mathcal J^{\vee}\otimes \mathcal H^{-n}, \mathcal V)\, ,\\
  \mathcal V^+\to \mathfrak M,\,\,\,\, & & \mathcal V^+=\text{Ext}^{\bullet}_{\,\widehat\pi}(\mathcal J^{\vee}, \mathcal V\otimes \mathcal J^{\vee})\, .
\end{eqnarray*}                                           
Pullbacks from the factors were suppressed in the expressions above. In particular, the above constructions make sense and will be used for
bundles $\mathcal V$ pulled back from $X$.

By relative duality, we have \begin{eqnarray}\label{rd} \mathcal V_n' & =& \text{Ext}^{\bullet}_{\,\widehat \pi} (\mathcal J^{\vee}\otimes \mathcal H^{-n}, \mathcal V)\\ &=& \text{Ext}^{\bullet}_{\,\widehat \pi} (\mathcal V, \mathcal
                                                                                                                                                                              J^{\vee}\otimes   \mathcal H^{-n}\otimes \omega_{\mathcal C/\BB})^{\vee}[1]\nonumber\\ &=&{\mathbf R} \widehat \pi_{\star} \left(\mathcal J^{\vee}\otimes \mathcal H^{-n}\otimes \mathcal V^{\vee} \otimes \omega_{\mathcal C/\BB}\right)^{\vee}[1]\nonumber \\ &=&\left(\mathcal V^{\vee}\otimes \omega_{\mathcal C/\BB}\right)_n^{\vee}[1]\, .
\nonumber\end{eqnarray}
The above constructions make sense for $K$-theory classes $\mathcal V$ as well.
\end{itemize}

Returning to Lemma \ref{l1}, we now
compute the pullbacks of the various tautological structures under the morphism $$j: (\mathcal C/\BB)^{[n]}\to X^{[n]}\, .$$

\begin{lemma} \label{l2}There are $K$-theory classes $\alpha, \beta$ on $\mathcal C$ and $\gamma$ on $\mathfrak M$ for which
$$j^{\star}\left(TX^{[n]}-\left(M^{[n]}\right)^{\vee}\right)=\gamma+\alpha_n\cdot \zeta_n^{-1}+\left(\beta_n\right)^{\vee}\cdot \zeta_n$$ over $(\mathcal C/\mathcal B)^{[n]}\to \mathbb P_n\to \mathfrak M$. Furthermore, $\alpha$ has rank $-1$ and $\beta$ has rank $0$. 
\end{lemma}
\proof We compute the two pullbacks separately. 
\begin{itemize} 
\item [(i)] First, recall  $$M^{[n]}={\mathbf R}\text{pr}_{\star} (M\otimes \mathcal O_{\mathcal W_n})$$ where $\mathcal W_n$ denotes the universal subscheme on $X^{[n]}\times X$ and
  $$ \text{pr}: X^{[n]} \times X \rightarrow X^{[n]}\, .$$
  The pullbacks on $M$ are omitted. 

  The pullback under $j$ is computed via the fibers of $$\pi:(\mathcal C/\BB)^{[n]}\times_{\BB}
  \mathcal C\to (\mathcal C/\BB)^{[n]}\, .$$
  We find $$j^{\star} M^{[n]}={\mathbf R}\pi_{\star} (M\otimes \mathcal O_{\mathcal Z_n})\, .$$
  Writing in $K$-theory $$\mathcal O_{\mathcal Z_n}=\mathcal O - \mathcal I_{\mathcal Z_n}=\mathcal O -\mathcal J^{\vee} \cdot \mathcal H^{-n} \cdot \zeta_n^{-1}$$ via equation \eqref{ideal}, we obtain 
  \begin{equation}\label{comput}j^{\star} M^{[n]}=\overline M - M_n \cdot \zeta_n^{-1}\, .\end{equation}
  Here, we have used the notations introduced in (A) and (B) above
  applied to the line bundle $M\to \mathcal C\to X$.

  \vspace{5pt}
\item [(ii)] We now turn to $j^{\star}TX^{[n]}.$ The alternating sum $$\mathcal O^{[n]}- TX^{[n]}+\left((K_X)^{[n]}\right)^{\vee}$$ computes fiber by fiber the complex $$\text{Ext}^0(\mathcal O_W, \mathcal O_W)-\text{Ext}^1(\mathcal O_W, \mathcal O_W)+\text{Ext}^2(\mathcal O_W, \mathcal O_W)$$ for subschemes $W$ of $X$. In families, 
  \begin{equation}\label{pba}j^{\star}\left(\mathcal O^{[n]}- TX^{[n]}+\left((K_X)^{[n]}\right)^{\vee}\right)=
    j^{\star} \Ext_{X}^{\bullet}(\mathcal O_{\mathcal W_n}, \mathcal O_{\mathcal W_n})\end{equation} where the subscript $X$ indicates the relative Ext's over the projection $$\text{pr}: X^{[n]}\times X\to X^{[n]}.$$ 

We seek to relate the relative $\Ext_{X}^{\bullet}$ and $\Ext_{\mathcal C/\BB}^{\bullet}$ where the second Ext is computed via the projection $$\pi:(\mathcal C/\BB)^{[n]}\times_{\BB}\mathcal C\to (\mathcal C/\BB)^{[n]}\, .$$
The key identity is \begin{equation}\label{rt}j^{\star} \Ext_{X}^{\bullet}(\mathcal O_{\mathcal W_n}, \mathcal O_{\mathcal W_n})=\Ext_{\mathcal C/\BB}^{\bullet}(\mathcal O_{\mathcal Z_n}, \mathcal O_{\mathcal Z_n}) -
  \Ext_{\mathcal C/\BB}^{\bullet}(\mathcal O_{\mathcal Z_n}\otimes \mathcal N, \mathcal O_{\mathcal Z_n})\,
  .\end{equation}  
Here, $\mathcal N$ is the conormal bundle of the inclusion $$\iota: \mathcal C\hookrightarrow \mathcal B \times X,$$ so that $$0\to \mathcal N\to \iota^{\star}\Omega_{\mathcal B\times X/\mathcal B} \to \Omega_{\mathcal C/\mathcal B}\to 0$$ or equivalently $$0\to \mathcal N\to j^{\star}\Omega_X\to \Omega_{\mathcal C/\mathcal B}\to 0\, .$$ 
Equation \eqref{rt} is the relative analogue of \cite[Lemma 3.42]{T} which gives the exactness of the sequence $$\text{Ext}^i_{\mathcal C_b}(\mathcal O_\xi, \mathcal O_\xi)\to \text{Ext}^i_{X}(\mathcal O_\xi, \mathcal O_\xi)\to \text{Ext}^{i-1}_{\mathcal C_b}(\mathcal O_\xi \otimes \mathcal N, \mathcal O_\xi)\to \text{Ext}^{i+1}_{\mathcal C_b}(\mathcal O_\xi, \mathcal O_\xi)\to \ldots$$ for subschemes $\xi\subset \mathcal C_b.$ To apply \cite{T}, we observe that $\mathcal N\big\rvert_{\mathcal C_b}$ is the conormal bundle of $\mathcal C_b\hookrightarrow X$,
which follows by restricting the defining exact sequence to $\mathcal C_b$.\footnote {We have $$0\to \mathcal Tor^{1}_{\mathcal C\,} (\Omega_{\mathcal C/\mathcal B}, \mathcal O_{\mathcal C_b}) \to \mathcal N\big\rvert_{\mathcal C_b}\to \Omega_{X}\big\rvert_{\mathcal C_b}\to \Omega_{\mathcal C_b}\to 0\, .$$
  $\mathcal Tor^1$ is supported on the finitely many singularities of $\mathcal C_b$. Since $\mathcal N\big\rvert_{\mathcal C_b}$ is locally free,
  $\mathcal Tor^1$ vanishes. Therefore, $\mathcal N\big\rvert_{\mathcal C_b}$ is the conormal bundle.}
\vskip.1in
With \eqref{rt} understood, and by invoking \eqref{pba}, we find \begin{equation} \label{pbaa} j^{\star} TX^{[n]} = j^{\star} \mathcal O^{[n]} + j^{\star} \left((K_X)^{[n]}\right)^{\vee} -
  \Ext_{\mathcal C/\BB}^{\bullet}(\mathcal O_{\mathcal Z_n}, \mathcal O_{\mathcal Z_n}) + \Ext_{\mathcal C/\BB}^{\bullet}(\mathcal O_{\mathcal Z_n}\otimes \mathcal N, \mathcal O_{\mathcal Z_n})\, .\end{equation}
The calculations in (i), specifically \eqref{comput}, yield the first two terms above $$j^{\star}  \mathcal O^{[n]} = \overline {\mathcal O} - \mathcal O_n \cdot \zeta_n^{-1}$$ and $$j^{\star} \left((K_X)^{[n]}\right)^{\vee}= {\overline {K_X}}^{\vee}-\left(\left(K_X\right)_n\right)^{\vee} \cdot\zeta_n\, .$$ 
We examine the last two terms in \eqref{pbaa}.
Substituting $$\mathcal O_{\mathcal Z_n}=\mathcal O - \mathcal J^{\vee} \mathcal H^{-n}\zeta_n^{-1}$$
yields \begin{eqnarray*} \Ext_{\mathcal C/\BB}^{\bullet}(\mathcal O_{\mathcal Z_n}, \mathcal O_{\mathcal Z_n})&=& \Ext_{\mathcal C/\BB}^{\bullet}(\mathcal O - \mathcal J^{\vee} \mathcal H^{-n}\zeta_n^{-1}, \mathcal O - \mathcal J^{\vee} \mathcal H^{-n}\zeta_n^{-1})\\ &=& \overline{\mathcal O} + \Ext_{\mathcal C/\BB} ^{\bullet}(\mathcal J^\vee, \mathcal J^\vee) - {\bf R}\pi_{\star} (\mathcal J^{\vee} \otimes \mathcal H^{-n})\cdot \zeta_n^{-1}\\
         & &  \ \ \ -\Ext^{\bullet}_{\mathcal C/\BB} (\mathcal J^{\vee}\otimes \mathcal H^{-n}, \mathcal O) \cdot \zeta_n\\ &=&  \overline{\mathcal O} + \mathcal O^+- \mathcal O_n\cdot \zeta_n^{-1} -\mathcal O_n'\cdot \zeta_n\, . \end{eqnarray*}
       An entirely similar calculation shows that \begin{eqnarray*}
 \Ext_{\mathcal C/\BB}^{\bullet}(\mathcal O_{\mathcal Z_n}\otimes \mathcal N, \mathcal O_{\mathcal Z_n})=\overline{\mathcal N^{\vee}} + \left(\mathcal N^{\vee}\right)^+- \left(\mathcal N^{\vee}\right)_n\cdot \zeta_n^{-1} -\left(\mathcal N^{\vee}\right)_n'\cdot \zeta_n\end{eqnarray*} Collecting the last four equations into \eqref{pbaa} we find
$$j^{\star} TX^{[n]}  = ({\overline {K_X}}^{\vee}+\overline{\mathcal N^{\vee}}+\left(\mathcal N^{\vee}\right)^+-\mathcal O^+) - \left(\mathcal N^{\vee}\right)_n\cdot \zeta_n^{-1}+\left(\mathcal O_n'-\left(\mathcal N^{\vee}\right)_n'-\left((K_X)_n\right)^{\vee}\right)\cdot
\zeta_n\, .$$ 
\end{itemize}
\vskip.1in
From (i) and (ii), we find that 
\begin{eqnarray*}
j^{\star} \left(TX^{[n]}-\left(M^{[n]}\right)^{\vee}\right)&=&({\overline {K_X}}^{\vee}+\overline{\mathcal N^{\vee}}  -\overline {M}^{\vee}+\left(\mathcal N^{\vee}\right)^+-\mathcal O^+) - \left(\mathcal N^{\vee}\right)_n \cdot \zeta_n^{-1}\\ && \hskip.5in+\left(\mathcal O_n'-\left(\mathcal N^{\vee}\right)_n'+ (M_n)^{\vee}-\left(\left(K_X\right)_n\right)^{\vee} \right)\cdot \zeta_n
.\end{eqnarray*} Using relative duality \eqref{rd} for the last terms, we rewrite the above answer as 
\begin{eqnarray*}({\overline {K_X}}^{\vee}&+&\overline{\mathcal N^{\vee}} -\overline {M}^{\vee}+\left(\mathcal N^{\vee}\right)^+-\mathcal O^+) - \left(\mathcal N^{\vee}\right)_n \cdot \zeta_n^{-1}\\
&+&  \left(-\left(\omega_{\mathcal C/\mathcal B}\right)_n^{\vee}+\left(\mathcal N\otimes \omega_{\mathcal C/\mathcal B}\right)_n^{\vee}+ (M_n)^{\vee}-\left(\left(K_X\right)_n\right)^{\vee} \right)\cdot \zeta_n\, ,
  \end{eqnarray*}
which
establishes Lemma \ref{l2}. \qed

\vspace{10pt}
We return now  to the proof of Lemma \ref{l1}.
First, we have  $$\mathbb P_n=\mathbb P(\epsilon'_n)$$ where $$\epsilon'_n={\bf R}\widehat\pi_{\star} (\mathcal J\otimes \mathcal H^{n})=\text{Ext}^{\bullet}_{\widehat{\pi}}
(\mathcal J^{\vee}\otimes \mathcal H^{-n}, \mathcal O)=\mathcal O'_n$$ in the notation (B) above. By Lemma \ref{l2}, expression \eqref{e1} becomes \begin{equation}\label{e2}\int_{\mathbb P_n} c(\gamma+\alpha_n\cdot \zeta_n^{-1}+\left(\beta_n\right)^{\vee}\cdot \zeta_n)\end{equation} which we will  prove is of the form \eqref{special} for sufficiently large $n$.

The classes $\alpha$ and $\beta$ have ranks $-1$ and $0$ respectively.
Therefore,
$$\text{rank }\alpha_n=n+r_1, \,\,\text{rank }\beta_n=r_2,\,\,\,\text{rank } \epsilon'_n=n+r_3+1\, ,$$ for constants $r_1, r_2, r_3$. Let $m$ denote the dimension of $\mathfrak M$, and let $$d=m+r_3\, .$$ We obtain \begin{eqnarray*}\int_{\mathbb P_n} c(\gamma+\alpha_n\cdot \zeta_n^{-1}+\left(\beta_n\right)^{\vee}\cdot \zeta_n)&=&\sum_{u+v+w=n+d} \int_{\mathbb P_n} c_u(\gamma) \cdot c_v(\alpha_n\cdot \zeta_n^{-1}) \cdot c_w(\left(\beta_n\right)^{\vee}\cdot \zeta_n)\, .\end{eqnarray*} The usual formulas give \begin{eqnarray*} && c_v(\alpha_n\cdot \zeta_n^{-1})=\sum_{i=0}^{v}\binom{\text{rank } \alpha_n-i}{v-i} \cdot c_{i}(\alpha_n)\cdot c_1(\zeta_n)^{v-i}\cdot (-1)^{v-i}\\ &&c_w(\left(\beta_n\right)^{\vee}\cdot \zeta_n)=\sum_{j=0}^{w}\binom{\text{rank } \beta_n-j}{w-j}\cdot c_{j} (\left(\beta_n\right)^{\vee})\cdot c_1(\zeta_n)^{w-j}\, .\end{eqnarray*} We therefore are led to the expressions 
$$\sum_{u, v, w, i, j} (-1)^{v-i} \binom{n+r_1-i}{v-i} \binom{r_2-j}{w-j}\int_{\mathbb P_n} c_u(\gamma) \cdot c_{i}(\alpha_n) \cdot c_{j} (\left(\beta_n\right)^{\vee}) \cdot c_1(\zeta_n)^{v+w-i-j}\, .$$
Integrating out $\zeta_n$ over the fibers of $$\mathbb P_n\to \mathfrak M\, ,$$ we rewrite expressions \eqref{e1} and \eqref{e2} as \begin{equation}\label{e3}\sum_{u+i+j+k=m}  \sigma_{ijk}^{(n)}\int_{\mathfrak M} c_u(\gamma)\cdot c_i(\alpha_n) \cdot c_j (\left(\beta_n\right)^{\vee}) \cdot c_k(-\epsilon'_n)\, ,\end{equation} where $$\sigma_{ijk}^{(n)}=\sum_{v+w=i+j+k+(n+r_3)} (-1)^{v-i}\binom{n+r_1-i}{v-i} \binom{r_2-j}{w-j}\, .$$ The number of terms in this binomial sum could potentially grow with $n$. However, $i, j, k$ are bounded independently of $n$.  \vskip.1in

\begin{lemma} \label{l3} For all $i, j, k$, $\sigma_{ijk}^{(n)}$ is of the form \eqref{special}. 
\end{lemma}

\begin{lemma}\label{l4} There exists $M>0$ and $K$-classes $\mu^{(i)}_0, \ldots, \mu^{(i)}_M$ on $\mathfrak M$, for $1\leq i\leq 3$, satisfying
\begin{eqnarray*}
  \alpha_n&=&\sum_{\ell=0}^{M} n^\ell \mu_\ell^{(1)}\, ,\\
  \left(\beta_n\right)^{\vee}&=&\sum_{\ell=0}^{M} n^\ell \mu_\ell^{(2)}\, ,\\
  -\epsilon'_n&=&\sum_{\ell=0}^{M} n^\ell \mu_\ell^{(3)}\, .
                  \end{eqnarray*}
\end{lemma}

Assuming Lemmas \ref{l3} and \ref{l4}, we return to \eqref{e3} and the proof of Lemma \ref{l1}. For any classes $\mu_\ell$ on $\mathfrak M$ with $0\leq \ell\leq M$, we have
$$c_i\left(\sum_{\ell=0}^{M} n^\ell \mu_\ell\right)=\left[\prod_{\ell=0}^{M} \left(1+c_1(\mu_\ell)+c_2(\mu_{\ell})+\ldots\right)^{n^\ell}\right]_{(i)}\, .$$
Furthermore, 
$$\left(1+c_1(\mu_\ell)+c_2(\mu_{\ell})+\ldots\right)^{n^\ell}=\sum_{I} c_I(\mu_\ell) \binom{n^{\ell}}{I}\, ,$$ where $I$ is a multi-index.
Therefore $$c_i\left(\sum_{\ell=0}^{M} n^\ell \mu_\ell\right)=\sum c_{I_0}(\mu_0) c_{I_1} (\mu_1)\cdots c_{I_{M}} (\mu_{M}) \cdot \binom{1}{I_0}\binom{n}{I_1}\cdots \binom{n^{M}}{I_M}\, ,$$
where, for degree reasons, $$|I_0|+\ldots+|I_{M}|=i\, .$$ Thus $I_j$'s have entries bounded by $i\leq m$, and the above expression is therefore polynomial in $n$. As a result, \eqref{e3} becomes 
$$\sum_{u+i+j+k=m} \sigma_{ijk}^{(n)} \int_{\mathfrak M} c_u(\gamma) \cdot c_{i}\left(\sum_{\ell=0}^{M} n^\ell \mu_\ell^{(1)}\right)\cdot c_{j}\left(\sum_{\ell=0}^{M} n^\ell \mu_\ell^{(2)}\right)\cdot c_{k}\left(\sum_{\ell=0}^{M} n^\ell \mu_\ell^{(3)}
\right)\, $$ which is of the form \eqref{special} by Lemma \ref{l3} and the above observations.
The proof of Lemma \ref{l1} will therefore be complete once Lemmas \ref{l3} and \ref{l4}
are proven.
\qed \vskip.1in

\subsubsection{Proof of Lemma \ref{l3}.}
The notation $$v^{\text{new}}=v-i\, , \ \,w^{\text{new}}=w-j\, ,
\ \, a=r_1-i\, ,\ \, b=r_2-j\, , \ \,c=r_3+k\, $$
will be more convenient for us.
With the new conventions, the expression in Lemma \ref{t3}
becomes $$\sigma^{(n)}=\sum_{v+w=n+c} (-1)^{v}\binom{n+a}{v}\binom{b}{w}=\text{Coeff}_{x^{n+c}} (1-x)^{n+a}\cdot(1+x)^{b}\, .$$
We rewrite the above as a residue $$\sigma^{(n)}\ =\
\text{Res}_{\,x=0\,} \frac{(1-x)^{n+a}\cdot (1+x)^b}{x^{n+c+1}}\,dx\, .$$
We change variables $$y=\frac{1-x}{x}\implies x=\frac{1}{y+1}$$
so that the differential form transforms to
$$\omega= - y^{n+a}\cdot (y+2)^{b} \cdot (y+1)^{e} \,dy\, ,$$
for some constant $e$.
Thus $$\sigma^{(n)}\, =\, \text{Res}_{\,y=\infty\,}\, \omega
\, =\, -\text{Res}_{\,y=-1\,}\omega-\text{Res}_{\,y=-2\,}\omega\, ,$$ via the
Residue Theorem. There are no poles for $\omega$ at $y=0$ for $n$ sufficiently large. 

The residues at $$y=-1,\,\, y=-2$$ correspond to the two terms of \eqref{special}. Indeed, for $y=-1$, we have
\begin{eqnarray*}
  \text{Res}_{\,y=-1\,} \,y^{n+a}\cdot (y+2)^{b} \cdot (y+1)^{e} \,dy
  &=&\text{Res}_{\,z=0\,} \, (z-1)^{n+a} \cdot (z+1)^b\cdot z^e \,dz \\
  & = &\text{Coeff}_{z^{-e-1}}\,(z-1)^{n+a} \cdot (z+1)^b\\
  &=&(-1)^{n} \sum_{v+w=-e-1} (-1)^{v+a} \binom{n+a}{v}\binom{b}{w}\, .
\end{eqnarray*}
 The latter sum is finite, hence manifestly polynomial in $n$. A similar calculation shows that the residue at $y=-2$ is of the form $$(-2)^n \cdot \text {polynomial in }n,$$ completing the proof. \qed
\vskip.1in

\subsubsection{Proof of Lemma \ref{l4}. }
We present the argument for $\alpha_n$. The proofs of the other two statements are the same. Consider the class
$$x=\mathcal H-1$$
viewed in the $K$-theory of $\mathcal C$. Since $\mathcal C$ is nonsingular
and projective, the Chern character gives an isomorphism $$\text{ch :}\, K(\mathcal C)\otimes \mathbb Q\to A^{\star}(\mathcal C)\otimes \mathbb Q\, .$$ Clearly $\text{ch}(x)\in A^{>0}(\mathcal C)$, hence $\text{ch}(x)^{M}=0$ for degree reasons, for some $M>0$. So $$x^{M}=0$$ in $K$-theory, hence
$$(\mathcal H-1)^{M}=0\, .$$ 
We conclude $$\sum_{i=0}^{M} (-1)^i \cdot \binom{M}{i}\cdot \mathcal H^i=0$$ and therefore  $$\sum_{i=0}^{M} (-1)^i \cdot \binom{M}{i}\cdot \mathcal H^{i-n}\cdot \mathcal J^{\vee} \alpha=0$$ in $K(\mathfrak M\times_B\mathcal C)\, .$ Pushing
forward via $\widehat \pi$ to $\mathfrak M$, we obtain $$\sum_{i=0}^{M} (-1)^i \cdot \binom{M}{i}\ \cdot \alpha_{n-i}=0\, .$$ This linear recursion in the $\alpha$'s can be solved explicitly. Note that the characteristic equation $$\sum_{i=0}^{M} (-1)^i \cdot \binom{M}{i}\ \cdot r^{-i}=0$$ has $M$ repeated roots all equal to $1$.  \qed
\vskip.1in

\subsubsection{The case $B\neq \emptyset$}\label{beee}
The last step of the proof of Proposition \ref{p1a} is to
treat the case when the sequence $B\neq \emptyset$.

The argument is similar to the $B=\emptyset$ case, but for completeness we indicate the main points. We wish to prove that $$Z_{X, M}[a, B]=\sum_{n=0}^{\infty} q^n\cdot \int_{X^{[n]}}  c_{n-a}\left(\left(M^{[n]}\right)^{\vee}\right) \cdot c(TX^{[n]}) \cdot \frac{P\left(\left(M^{[n]}\right)^{\vee},B\right)}{c\left(\left(M^{[n]}\right)^{\vee}\right)}$$ is rational in $q$. Following the reasoning in Section \ref{uni2} and using Lemma \ref{l2}, it suffices to show  $$\int_{\mathbb P_n} c(\gamma+\alpha_n\cdot \zeta_n^{-1}+\left(\beta_n\right)^{\vee}\cdot \zeta_n)\cdot P\left(j^{\star} \left(M^{[n]}\right)^{\vee}, B\right)$$ is of the form \eqref{special}, for $n$ large enough. This is analogous to Lemma \ref{l1}. 

By \eqref{comput}, we have $$P_b\left(j^{\star} \left(M^{[n]}\right)^{\vee}\right)=P_b (\overline M^{\vee}) + P_b \left(\left(- M\right)_n^{\vee}\cdot \zeta_n\right),$$ where we have extended the definition of $P_b$ given in \eqref{pbe} to $K$-theory by linearity. We multiply out the $P_b$'s for the values of $b$ determined by the sequence $B$. Since $\overline M$ is a $K$-theory class on $\mathfrak M$, we can combine terms of the form $P_b (\overline M^{\vee})$ and $c(\gamma)$ into a single cohomology class $\lambda$ over $\mathfrak M$. We are led to expressions of the form 
\begin{equation}\label{eee}\int_{\mathbb P_n} \lambda\cdot c(\alpha_n\cdot \zeta_n^{-1}+\left(\beta_n\right)^{\vee}\cdot \zeta_n)\cdot P_{b_1}\left((-M)_n^{\vee}\cdot \zeta_n\right)\cdots P_{b_m}\left((-M)_n^{\vee}\cdot \zeta_n\right),\end{equation} 
\noindent for nonnegative integers $b_1, \ldots, b_m$. 

Recall from Section \ref{L1L1L1} that
  $$(-M)_n=-{\bf R}\widehat \pi_{\star}(M\otimes \mathcal J^{\vee}\otimes \mathcal H^{-n})\, .$$ By inspecting the fiber degree, we see that $(-M)_n$ is represented by a vector bundle for $n$  sufficiently large, and we write $x_1, \ldots, x_{r}$ for the Chern roots. The rank $r$ depends on $n$ linearly. We have $$P_b\left((-M)_n^{\vee}\cdot \zeta_n\right)=\sum_{i=1}^{r} \frac{1}{(1-x_i+c_1(\zeta_n))^b}\, .$$
We expand $$\frac{1}{(1-x_i+c_1(\zeta_n))^{b}}=\frac{1}{(1-x_i)^b} \cdot \left(1+\frac{c_1(\zeta_n)}{1-x_i}\right)^{-b}=\sum_{\ell=0}^{\infty} \binom{-b}{\ell} \cdot \frac{c_1(\zeta_n)^\ell}{(1-x_i)^{\ell+b}}\, .$$ Thus $$P_b\left((-M)_n^{\vee}\cdot \zeta_n\right)=\sum_{\ell=0}^{\infty} \binom{-b}{\ell}\cdot c_1(\zeta_n)^\ell\cdot \mathsf p_{\ell+b}\left((-M)_n\right)$$ where the classes $\mathsf p$ on $\mathfrak M$ have terms of mixed degrees. In fact, the $\mathsf p$'s are series in the Chern classes of the argument whose coefficients are independent of $n$. The only exception is the constant term which is linear in $n$ being equal to the rank.

Expression \eqref{eee} becomes \begin{eqnarray*}\sum_{\ell_1, \cdots, \ell_m} \binom{-b_1}{\ell_1}\cdots \binom{-b_{m}} {\ell_m}\cdot \int_{\mathbb P_n} \lambda &\cdot& c(\alpha_n\cdot \zeta_n^{-1}+\left(\beta_n\right)^{\vee}\cdot \zeta_n)  \cdot  c_1(\zeta_n)^{\ell_1+\ldots+\ell_m}\\ &\cdot& \mathsf p_{\ell_1+b_1} \left((-M)_n\right)\cdots \mathsf p_{\ell_m+b_m} \left((-M)_n\right).\end{eqnarray*} To go further, we apply the same reasoning that led to equation \eqref{e3}. Accounting for the extra $\mathsf p$'s and their prefactors, the above expression becomes 
$$\sum_{i, j, k, \vec a}  \sigma_{i, j, k, \vec a}^{(n)} \int_{\mathfrak M} \lambda \cdot c_{i}(\alpha_n) \cdot c_{j} (\left(\beta_n\right)^{\vee})\cdot c_k(-\epsilon'_n) \cdot \mathsf p_{a_1} \left((-M)_n\right)\cdots \mathsf p_{a_m} \left((-M)_n\right)$$ where $$\ell_1+b_1=a_1\, , \ldots, \ell_m+b_{m}=a_m\, .$$ The prefactor here equals 
$$\sigma_{i, j, k, \vec a}^{(n)}=\sum_{v, w,  \vec b, \vec \ell} \binom{-b_1}{\ell_1}\cdots \binom{-b_{m}} {\ell_m}\cdot(-1)^{v-i} \binom{n+r_1-i}{v-i} \binom{r_2-j}{w-j}\, ,$$ where in the summation we have \begin{equation}\label{vwvw}v+w+|\ell|=i+j+k+(n+r_3)\, ,\,\, \,\, \vec \ell+\vec b=\vec a\, .\end{equation} Each integral over $\mathfrak M$ is polynomial in $n$. Indeed, dimension constraints select only finitely many homogeneous pieces from the $\mathsf p$'s and from $\lambda$, of bounded degree. We then argue by invoking Lemma \ref{l4} applied to $$\alpha_n\, ,\  \left(\beta_n\right)^{\vee} ,\ -\epsilon'_n\, ,\, (-M)_n$$ combined with the analysis that followed the statement of the Lemma. 

To conclude, it remains to prove $\sigma^{(n)}$ is of the form \eqref{special}, the analogue of Lemma \ref{l3}. We have \begin{eqnarray*}\sigma_{i, j, k, \vec a}^{(n)}&=&\sum_{v, w, \vec b, \vec \ell} \binom{-b_1}{\ell_1}\cdots \binom{-b_{m}} {\ell_m}\cdot(-1)^{v-i} \binom{n+r_1-i}{v-i} \binom{r_2-j}{w-j}\\ &=& \sum_{v, w, \vec b, \vec \ell} \binom{b_1+\ell_1-1}{\ell_1}\cdots \binom{b_{m}+\ell_m-1} {\ell_m}\cdot(-1)^{v-i+|\ell|} \binom{n+r_1-i}{v-i} \binom{r_2-j}{w-j}\\&=& \sum_{v, w, |\ell|} (-1)^{v-i+|\ell|}\binom{|a|-m}{|\ell|} \binom{n+r_1-i}{v-i}  \binom{r_2-j}{w-j} \end{eqnarray*} where the Vandermonde identity was used in the last line to sum over $\ell_1+\ldots+\ell_m=|\ell|$. Writing $|\ell|+v=v'$, and using the Vandermonde identity one more time, we obtain $$\sigma_{i, j, k, \vec a}^{(n)}=\sum_{v'+w=n+\text{const}} (-1)^{v'-i} \binom{(|a|-m)+(n+r_1-i)}{v'-i} \binom{r_2-j}{w-j}\, .$$ This is exactly the type of expression considered in Lemma \ref{l3}.
The proof of Proposition \ref{p1a} is complete. \qed

\subsubsection{Example} 
We illustrate the methods used in the proof of Proposition \ref{p1a} with the computation of the series \begin{equation}\label{zz}Z=\sum_{n=1}^{\infty} q^n \int_{X^{[n]}} c_{n-1}\left(\left(M^{[n]}\right)^{\vee}\right)\cdot
  c(TX^{[n]})\cdot s\left(\left(M^{[n]}\right)^{\vee}\right)\end{equation}
in the special case $$X=B\times F\to B \, \text { and }\,  M=\mathcal O_B(1)\, ,$$
where $B=F=\mathbb P^1$. The family of curves in the fiber class $\mathcal C\to B$ is isomorphic to the surface $X\to B$, and the relative Hilbert scheme of points is the product
$$\left(\mathcal C/B\right)^{[n]}=B\times \mathbb P^n\, .$$ By \eqref{ideal}, the universal subscheme $$\mathcal Z_n\hookrightarrow \mathcal C\times_{B}\left(\mathcal C/B\right)^{[n]}=X\times \mathbb P^n$$ satisfies $$I_{\mathcal Z_n} = \mathcal O_F(-n)\otimes  {\mathcal O}_{\mathbb P^n}(-1)\, .$$ 
We represent $$c_{n-1}\left(M^{[n]}\right)=\left(\mathcal C/B\right)^{[n]}= \left[B\times \mathbb P^n\right],$$ so that $$Z=\sum_{n=1}^{\infty} q^n (-1)^{n-1} \int_{B\times \mathbb P^n} j^{\star} c\left(TX^{[n]}-\left(M^{[n]}\right)^{\vee}\right).$$ Here, we continue to write $$j: B\times \mathbb P^n\to X^{[n]}$$ for the natural morphism. Let $$\pi: X\times \mathbb P^n \to B\times \mathbb P^n$$ denote the projection. 
We compute the tautological structures \begin{eqnarray*} j^{\star} M^{[n]}&=&{\mathbf R}\pi_{\star} (M \otimes \mathcal O_{\mathcal Z_n})\\&=&{\mathbf R}\pi_{\star} (M - M\otimes \mathcal O_F(-n)\otimes {\mathcal O}_{\mathbb P^n}(-1) )\\&=&\mathcal O_B(1) + \mathbb C^{n-1} \otimes \mathcal O_B(1)\otimes \mathcal {\mathcal O}_{\mathbb P^n}(-1)\, .\end{eqnarray*} Similarly, \begin{eqnarray*} j^{\star} \mathcal O_X^{[n]}&=&\mathcal O + \mathbb C^{n-1} \otimes {\mathcal O}_{\mathbb P^n}(-1)\, ,\\ j^{\star} \left(K_X\right)^{[n]}&=&-\,\mathcal O_B(-2) +\mathbb C^{n+1}\otimes \mathcal O_B(-2)\otimes {\mathcal O}_{\mathbb P^n}(-1)\, . 
                                                                                                                                                                                                                                                                                                                                                                                       \end{eqnarray*}By \eqref{pbaa}, we have $$j^{\star} TX^{[n]}=j^{\star} \mathcal O_X^{[n]}+  j^{\star}\left( \left(K_X\right)^{[n]}\right)^{\vee} - \text{Ext}^{\bullet}_{\mathcal C/B} (\mathcal O_{Z_n}, \mathcal O_{Z_n})+  \text{Ext}^{\bullet}_{\mathcal C/B} (\mathcal O_{Z_n}, \mathcal O_{Z_n})\otimes \mathcal O_B(2)\, .$$ Here, we have used
$$\mathcal N=\Omega_{B}=\mathcal O_B(-2)\, .$$ Furthermore,
\begin{eqnarray*}\text{Ext}^{\bullet}_{\mathcal C/B} (\mathcal O_{Z_n}, \mathcal O_{Z_n})&=&\text{Ext}^{\bullet}_{\mathcal C/B} (\mathcal O-\mathcal O_F(-n)\otimes {\mathcal O}_{\mathbb P^n}(-1), \mathcal O-\mathcal O_F(-n)\otimes {\mathcal O}_{\mathbb P^n}(-1))\\&=& \mathbb C^2\otimes \mathcal O_B +\mathbb C^{n-1}\otimes {\mathcal O}_{\mathbb P^n}(-1)- \mathbb C^{n+1}\otimes {\mathcal O}_{\mathbb P^n}(1)\, .\end{eqnarray*}
After substituting, we find \begin{eqnarray*} j^{\star}  \left(TX^{[n]}-\left(M^{[n]}\right)^{\vee}\right)= \left(-\mathcal O_B  - \mathcal O_B(-1)+\mathcal O_B(2)\right) + \mathbb C^{n-1} \otimes \mathcal O_B(2) \otimes {\mathcal O}_{\mathbb P^n}(-1) \\ + \, \mathbb C^{n+1} \otimes {\mathcal O}_{\mathbb P^n}(1)-\mathbb C^{n-1}\otimes \mathcal O_B(-1)\otimes {\mathcal O}_{\mathbb P^n}(1)\, .
                            \end{eqnarray*} With $h, \zeta$ denoting the hyperplane classes on $B$ and $\mathbb P^n$, we arrive at the integral $$\int_{B\times \mathbb P^n} \frac{1+2h}{1-h}\cdot (1+\zeta)^{n+1}\cdot \left(\frac{1-\zeta+2h}{1+\zeta-h}\right)^{n-1}.$$ The last expression equals $$\text{Coeff}_{\,h\zeta^n\,} \ \frac{1+2h}{1-h}\cdot (1+\zeta)^{n+1}\cdot \left(\frac{1-\zeta+2h}{1+\zeta-h}\right)^{n-1}= (-1)^{n} \cdot (4n-10)\, .$$
Hence, we can write \eqref{zz} as $$Z=\sum_{n=1}^{\infty} q^n (10-4n)=\frac{q(6-10q)}{(1-q)^2}\, .$$

For another example, if $X$ is a $K3$ surface and $M=\mathcal O_X,$ the series $$Z=\sum_{n=1}^{\infty} q^n \int_{X^{[n]}} c_{n-1}\left(\left(M^{[n]}\right)^{\vee}\right)\cdot c(TX^{[n]})\cdot s\left(\left(M^{[n]}\right)^{\vee}\right)=\frac{24q}{(1-q)^2}$$ was computed in Proposition 40 of \cite {OP1}. 

Evaluating \eqref{zz} in closed form for all pairs $(X, M)$ is likely possible. 

\section{Descendent series of punctual Quot schemes:
  Theorem \ref{t2}} \label{sss3}

\subsection{Overview} The goal here is to prove Theorem \ref{t2}. Throughout  Section \ref{sss3}, we set $\beta=0$. We will
establish the rationality of the descendent series  \begin{eqnarray*}Z_{X,\,N}(\alpha_1,\ldots, \alpha_\ell \,|\,k_1, \ldots, k_{\ell})=\sum_{n=0}^{\infty} q^n \cdot \int_{\left[\mathsf {Quot}_X(\mathbb C^N, n)\right]^{\text{vir}}}
  \mathsf{ch}_{k_1}(\alpha_1^{[n]}) \cdots \mathsf{ch}_{k_\ell}(\alpha_\ell^{[n]})\, c(T^{\text{vir}}\mathsf {Quot})\, .\end{eqnarray*}
Our argument follows the strategy of the proof of Theorem 18 of \cite {OP1}. 

\subsection {Proof of Theorem \ref{t2}} We will explain shortly that for fixed $\text{rank }\alpha_i=r_i$, the series $Z_{X,\,N}(\alpha_1,\ldots, \alpha_\ell \,|\,k_1, \ldots, k_{\ell})$ is given by universal expressions in the Chern classes of the $\alpha_i$'s. Furthermore, for each $k_1, \ldots, k_{\ell}$, the series $$Z_{X,\,N}(\alpha_1,\ldots, \alpha_\ell \,|\,k_1, \ldots, k_{\ell})$$ is additive in $\alpha_1, \ldots, \alpha_\ell$ separately. 
Thus, invoking the splitting principle, it suffices to assume that $$\text{rank }\alpha_i=1 \text{ for all }1\leq i\leq \ell.$$ The proof below can also be directly written for $\alpha_i$'s of arbitrary ranks, at the expense of more complicated notation.

Since the Chern character is polynomial in the Chern classes, we equivalently consider the series $$Z=\sum_{n=0}^{\infty} q^n \cdot \int_{\left[\mathsf {Quot}_X(\mathbb C^N, n)\right]^{\text{vir}}}
 c_{k_1}(\alpha_1^{[n]}) \cdots c_{k_\ell}(\alpha_\ell^{[n]})\, \cdot c(T^{\text{vir}}\mathsf {Quot})\, .$$ Let $x_1, \ldots, x_{\ell}$ be formal variables. Write $$c_x=1+x c_1+ x^2 c_2 + \ldots$$ for the total Chern class, and set $$W=\sum_{n=0}^{\infty} q^n \cdot \int_{\left[\mathsf {Quot}_X(\mathbb C^N, n)\right]^{\text{vir}}} c_{x_1}(\alpha_1^{[n]}) \dots c_{x_\ell}(\alpha_\ell^{[n]})\cdot c(T^{\text{vir}}\mathsf {Quot})\, .$$ The series $Z$ is found by extracting the coefficient of $x_1^{k_1}\cdots x_{\ell}^{k_{\ell}}$ in $W$: $$Z=\frac{1}{k_1!}\cdots \frac{1}{k_\ell!}\cdot \frac{\partial^{k_1}}{\partial^{k_1}x_1} \ldots \frac{\partial^{k_\ell}}{\partial^{k_\ell}x_\ell} W\bigg\rvert_{x_1=\ldots=x_{\ell}=0}.$$ 
As in \cite{OP1}, we have a factorization
$$W = {\mathsf A}^{K_X^2} \cdot {\mathsf B}^{\chi(X)} \cdot \prod_{i=1}^{\ell} {\mathsf C}_i^{c_1(\alpha_i).K_X} \cdot {\mathsf D}_i^{c_1(\alpha_i)^2} \cdot \mathsf E_i^{c_2(\alpha_i)}\cdot \prod_{1\leq i<j\leq \ell} \mathsf F_{ij}^{c_1(\alpha_i)\cdot c_1(\alpha_j)}$$
for universal series $\mathsf A, \mathsf B, \mathsf C_i, \mathsf D_i, \mathsf E_i, \mathsf F_{ij}$ that depend on $q$ and $x_m$. We study the rationality of these series and of their $x_m$-derivatives. 

To this end, we pick convenient geometries. Take a nonsingular projective surface $X$ which admits a nonsingular
connected canonical curve $$\iota:C \hookrightarrow
 X$$ of genus $g$. We move the calculation to the punctual Quot scheme of the curve $C$: 
 $$\iota: \mathsf {Quot}_{C} (\mathbb C^N, n)\to \mathsf{Quot}_{X} (\mathbb C^N, n)\, .$$
By  \cite[Lemma 34]{OP1}, we have $$\iota_{\star} \left[\mathsf {Quot}_{C} (\mathbb C^N, n)\right]=(-1)^n  \left[\mathsf{Quot}_{X} (\mathbb C^N, n)\right]^{\mathrm{vir}}.$$ Furthermore, as remarked in equation (42) of \cite {OP1}, in $K$-theory we have the decomposition \begin{equation}\label{splt}\iota^{\star} T^{\mathrm{vir}} \text{Quot}_X(\mathbb C^N, n)=T \text{Quot}_C(\mathbb C^N, n)+
  \mathcal T_n\, .\end{equation}
Here, $\mathcal T_n\to \text{Quot}_C(\mathbb C^N, n)$ is the virtual bundle given pointwise over the quotient $$\mathbb C^N\otimes \mathcal O_C\to Q$$ 
by the expression $$\mathcal T_n=\text{Ext}^{\bullet}_C(Q, Q\otimes \Theta)\, ,$$ 
where $\Theta=N_{C/X}$ is the associated theta characteristic. 
As a consequence, we have 
 $$W=\sum_{n=0}^{\infty} q^n (-1)^n\cdot \int_{\mathsf {Quot}_C({\mathbb C}^N,n)} c_{x_1}(\iota^*\alpha_1) \cdots c_{x_\ell}( \iota^* \alpha_n) \cdot c(T\mathsf{Quot}_C) \cdot c(\mathcal T_n)\, .$$
The above expression does not depend on the surface $X$, which we will ignore from now on. It follows then that $$\mathsf B=1\, ,\,\,\, \mathsf {D}_i=1\, ,\,\,\, \mathsf E_i=1,\,\,\, \mathsf F_{ij}=1.$$ Therefore, for $\beta_i=\iota^{\star} \alpha_i$, we have 
$$W=\sum_{n=0}^{\infty} q^n (-1)^n\cdot \int_{\mathsf {Quot}_C({\mathbb C}^N,n)} c_{x_1}(\beta_1^{[n]}) \cdots c_{x_\ell}( \beta_\ell^{[n]}) \cdot c(T\mathsf{Quot}_C) \cdot c(\mathcal T_n)$$ with the factorization 
 $$W=\mathsf A^{g-1} \cdot \mathsf C_1^{\,c_1(\beta_1)} \cdots \mathsf C_{\ell}^{\,c_1(\beta_\ell)}\, .$$ 
 
We will establish that the $x_m$-derivatives of the series $\mathsf A$ and $\mathsf C_m$ are rational in $q$ after setting the $x$'s to $0$. To study these series, we may pick again convenient geometries: $$C=\mathbb P^1,\,\,\,\beta_i=\mathcal O_{\mathbb P^1}(d_i)$$ for arbitrary integers $d_i$. Therefore \begin{equation}\label{ys}
W=\sum_{n=0}^\infty q^n (-1)^n \cdot \int_{\mathsf {Quot}_{{\mathbb P}^1}(\mathbb C^N,n)} c_{x_1}(\OO(d_1)^{[n]})  \cdots c_{x_\ell} (\OO(d_\ell)^{[n]}) \cdot c(T\mathsf{Quot}_{\mathbb P^1})  \cdot c(\mathcal T_n)\, .\end{equation} It suffices to show the rationality of the $x_m$-derivatives of $W$. 

We use Atiyah-Bott equivariant localization to compute \eqref{ys}. We let the torus $\mathbb C^{\star}$ act on $\mathbb C^N$ with weights $$w_1\, , \ldots, w_N\, ,$$ thus inducing an action on $\mathsf {Quot}_{{\mathbb P}^1}(\mathbb C^N,n)$. The fixed loci were noted in \cite {OP1} to be isomorphic to $$C^{[n_1]}\times \cdots \times C^{[n_N]}=\mathbb P^{n_1}\times \cdots \times \mathbb P^{n_N}$$ for partitions $$n_1+\ldots+n_N=n\, .$$ Equivariant localization applied to \eqref{ys} thus yields 
\begin{equation}\label{yss}W=\sum_{n=0}^{\infty} q^n \cdot \sum_{n_1+\ldots+n_N=n} \int_{\mathbb P^{n_1}\times\cdots \times \mathbb P^{n_N}}
  \mathsf {Contr}(n_1, \ldots, n_N)\, .\end{equation} The expression $$\mathsf {Contr}(n_1, \ldots, n_N)$$ encodes the contribution of the fixed loci. In the absence of the descendent classes $c_x(\OO(d)^{[n]})$, the contribution was determined explicitly in \cite {OP1} in the proof of Theorem 18, via a calculation of the normal bundles of the fixed loci. The answer is \begin{equation*}
\mathsf {Contr}(n_1, \ldots,n_N)=(-1)^{nN+\binom{N}{2}} \cdot \Phi_1(h_1)^{n_1} \cdots \Phi_N(h_N)^{n_N} \cdot \Psi(h_1, \ldots, h_n)\, 
\end{equation*}
for the rational functions 
$$
\Phi_i(h_i)=\prod_{j=1}^{N} (1-h_i+w_i-w_j)\cdot \prod_{j\neq i} (h_i+w_j-w_i)^{-1}\, ,$$ 
\begin{align*}\Psi= \prod_{i<j} (h_i-h_j+w_j-w_i)^2 
  \cdot \prod_{i,j} (1+h_i+w_j-w_i)\cdot (1+h_i-h_j+w_j-w_i)^{-1} \\
  \cdot  \prod_{j\neq i}(h_i+w_j-w_i)^{-1}\, .
\end{align*} We must modify these
rational functions to account for the descendent insertions. 

We will use Lemma 27 of \cite {OP1}. For $N=1$, over $C^{[n]}=\mathbb P^n$, the tautological classes can be expressed in $K$-theory as $$\OO(d)^{[n]}=(d+1)\cdot  \mathcal O_{\mathbb P^n} + (-d+n-1)\cdot \mathcal O_{\mathbb P^n}(-1)\, .$$ Letting the torus act on $\OO(d)^{[n]}$ fiberwise with weight $w$ and writing $h$ for the hyperplane class, we obtain
$$ c_x(\OO(d)^{[n]}[w]) = (1+xw)^{d+1} \cdot (1+x(w - h))^{-d+n-1}.$$ For $N>1$, restricting  $\OO(d)^{[n]}$ to the fixed locus $$\mathbb P^{n_1}\times \cdots \times \mathbb P^{n_N}$$ yields $$c_x(\OO(d)^{[n]})=\prod_{i=1}^{N} (1+xw_i)^{d+1} \cdot \prod_{i=1}^{N} (1+x(w_i-h_i))^{-d+n_i-1},$$ with $h_i$ denoting the hyperplane classes of each factor. The new contributions to \eqref{yss} thus become \begin{equation*}
\mathsf {Contr}(n_1, \ldots,n_N)=(-1)^{nN+\binom{N}{2}} \cdot \Phi_1(h_1)^{n_1} \cdots \Phi_N(h_N)^{n_N} \cdot \Psi(h_1, \ldots, h_n)\, 
\end{equation*}
for the new rational functions 
$$
\Phi_i(h_i)=\prod_{j=1}^{N} (1-h_i+w_i-w_j)\cdot \prod_{j\neq i} (h_i+w_j-w_i)^{-1}\, \cdot  \prod_{m=1}^\ell (1 + x_m(w_i - h_i))\, ,$$ 
\begin{align*}
\Psi = &\prod_{i<j} (h_i - h_j +w_j - w_i)^2 \cdot \prod_{i,j} (1 + h_i + w_j - w_i) \cdot (1 + h_i - h_j + w_j -w_i)^{-1}  \\
&\cdot \prod_{j \neq i} (h_i + w_j - w_i)^{-1} \cdot 
\prod_{i} \prod_{m=1}^{\ell} (1+x_m(w_i-h_i))^{-1-d_m} \cdot (1+x_m w_i)^{d_m+1}\, .
\end{align*}
We conclude
\begin{multline*}
  W=\sum_{n=0}^{\infty} q^n\,  (-1)^{nN+\binom{N}{2}}\\
  \cdot \sum_{n_1+\ldots+n_N=n} \left[h_1^{n_1}\cdots h_N^{n_N}\right]\, \Phi_1(h_1)^{n_1} \cdots \Phi_N(h_N)^{n_N} \cdot \Psi(h_1, \ldots, h_n)\, .
  \end{multline*} The brackets in the above
series  are used to denote the coefficient of the relevant monomial. 

By Lagrange-B\"urmann's formula \cite {G}, we obtain $$W=(-1)^{\binom{N}{2}}\cdot \frac{\Psi}{K}\left(h_1, \dots, h_N\right)$$
where as usual 
$$K = \prod_{i=1}^N \left(1-h_i\cdot \frac{\Phi_i'(h_i)}{\Phi_i(h_i)}\right),$$
and $h_i$ is the solution to the equation
$$(-1)^N q= \frac{h_i}{\Phi_i(h_i)}$$
satisfying $h_i(q=0)=0$. At the end, we also set the equivariant weights $w_i$ equal to zero.

We define 
$$\mathsf X(g) = \prod_{j=1}^N \frac{g+w_j}{1-g - w_j} \cdot \prod_{m=1}^\ell \frac{1}{1-x_mg}\, .$$ We then have 
$$\frac{h_i}{\Phi_i(h_i)} = \mathsf X(h_i-w_i)\, .$$ Let $g_1, \dots, g_N$ be the solutions{\footnote{There are other roots  which we will deal with later.
     See
    equation \eqref{eq:P-fac}.
  }}
to $$ \mathsf X(g)=(-1)^{N} q$$ with initial conditions $$g_i(q=0) = -w_i\, .$$ 
The $g_i$ are power series in $q$ whose coefficients are rational functions in $\{w_i\}$ and $\{x_m\}$.  Thus $$h_i=g_i + w_i$$ is a solution to $$\frac{h_i}{\Phi_i(h_i)} = (-1)^{N} q$$ with $h_i(q=0)=0$. 

We can easily check 
\begin{align*}
    K(h_1, \dots, h_N) \, = \, \prod_{i=1}^{N} h_i \frac{d}{dh_i} \log \frac{h_i}{\Phi_i(h_i)} 
  \,  =\, \prod_{i=1}^{N} h_i \frac{d}{dg_i} \log \mathsf X(g_i)\, .
\end{align*}
Furthermore, 
\begin{align*}
\Psi(h_1, \dots, h_N) &= &\prod_{i=1}^{N} h_i\cdot \prod_{i<j} (g_i-g_j)^2 \cdot \prod_{i,j} (1+g_i + w_j) \cdot (1+g_i - g_j)^{-1}\cdot (g_i +w_j)^{-1} \\
&& \cdot \prod_{i=1}^N \prod_{m=1}^{\ell} (1-g_ix_m)^{-1-d_m}\cdot (1+x_mw_i)^{d_m+1}\, .
\end{align*}
The expressions for $K$ and $\Psi$ are evidently symmetric in $g_1, \ldots, g_N$, except for the factor $\prod_{i=1}^{N} h_i$ which appears in both. Hence the quotient $$\frac{\Psi}{K}(h_1, \dots, h_N)$$ can be expressed as a rational function in $\{g_i\}$, $\{w_i\}$,  and $\{x_m\}$ which is symmetric in the $\{g_i\}$. 

\vskip.1in

We rewrite the equation $\mathsf X(g)=(-1)^{N} q$ as $\mathsf P(g)=0$, where
\begin{align}
 \mathsf P (g) &= \prod_{i=1}^{N}(g + w_i) - q \prod_{i=1}^{N} (g+w_i-1) \prod_{m=1}^{\ell} (1- x_m g) \nonumber \\
 			 &= \sum_{j=0}^{N+\ell} \mathsf P_j g^j\, . \label{eq:P-def}
\end{align}
The $\{g_i\}$ are roots of $\mathsf P$.  Hence, $\mathsf P$ factors as
\begin{align}
  \mathsf P (g) &= \prod_{i=1}^N \left(g-g_i\right) \cdot (f_\ell g^\ell + f_{\ell-1}g^{\ell-1} + \dots + f_0)  \nonumber \\
  &=(g^N + e_{N-1}g^{N-1} + \dots + e_0)\cdot (f_\ell g^\ell + f_{\ell-1}g^{\ell-1} + \dots + f_0) \label{eq:P-fac}
\end{align}
where $e_i$ is $(-1)^i$ times the $(N-i)$th elementary symmetric function in $g_i$, and the $f_m$ are  power series in $q$ with coefficients given by
rational functions of $\{w_i\}$ and $\{x_m\}$.  Now, setting $q=0$, we see that 
\begin{align*}
  \prod_{i=1}^{N} (g+w_i) &= \mathsf P(g, q=0) \\
  &= \prod_{i=1}^{N} (g - g_i(0)) \cdot (f_\ell(0) g^\ell + f_{\ell-1}(0)g^{\ell-1} + \dots + f_0(0)) \\
  &=\prod_{i=1}^{N} (g +w_i) \cdot (f_\ell(0) g^\ell + f_{\ell-1}(0)g^{\ell-1} + \dots + f_0(0))\, .
\end{align*}
It follows then that $f_m(0) = 0$ for $m>0$, and $f_0(0) = 1$.

We claim that both $\{e_i\}$ and $\{f_m\}$ are series in $q$ whose coefficients are polynomials in $\{w_i\}$ and $\{x_m\}$. We will abbreviate this by saying that they ``are polynomial". So far, we can see that this is true up to order 0. Let us assume, by induction, that this is true to order $p$.

For $m = \ell, \ell-1, \dots, 0$, we compare the coefficient of $g^{m+N}$ in
the expressions \eqref{eq:P-def} and \eqref{eq:P-fac}. We have
$$ f_{m} + \sum_{k\geq 1} f_{m+k}e_{N-k} = \mathsf P_{m+N}.$$
Here, $k$ goes up to the minimum of $N$ and $\ell-m$, and the sum is empty for $m = \ell$.
By inducting on $m$, we may assume that all the $f_{m+k}$ are polynomial to order $p+1$. Since these $f_{m+k}$ also have no constant term, the $(p+1)$st term of $e_{N-k}$ is not needed to compute the $(p+1)$st term of $f_{m+k}e_{N-k}$. In addition, $\mathsf P_{m+N}$ is known exactly and is polynomial. Hence, we see that $f_m$ is polynomial to order $p+1$.

Now, for $i=0, \dots, N-1$, we compare the  coefficients of $g^i$ in the expressions \eqref{eq:P-def} and \eqref{eq:P-fac}. We have
$$ e_if_0 + \sum_{k\geq 1} e_{i-k} f_k = \mathsf P_i.$$
By inducting on $i$, we may assume that the $e_{i-k}$ are polynomial to order $p+1$. We also know the $f_k$ are polynomial to order $p+1$. We know that $f_0$ starts with $1$, so it has a multiplicative inverse which is also polynomial to order $p+1$. It follows that $e_i$ is polynomial to order $p+1$. Our induction on $p$ is complete.
\vskip.1in

 We would like to see that any order derivative $$\frac{\partial^{k_1}}{\partial^{k_1}x_1} \ldots \frac{\partial^{k_\ell}}{\partial^{k_\ell}x_\ell}\bigg\rvert_{x_1=\ldots=x_{\ell}=0}$$
of $\Psi/K$ is a rational function in $q$, after setting the $w's$ to zero. This will follow from the observations below. 

Fix a rational function $\mathsf R$ of $\{g_i\}$, $\{w_i\}$, $\{x_m\}$ and $q$, which is symmetric in the $\{g_i\}$. Of course $\mathsf R$ can be rewritten as a rational function of $\{e_i\}$, $\{w_i\}$, $\{x_m\}$ and $q$. Setting $\{w_i=0\}$ and $\{x_m=0\}$ in \eqref{eq:P-def} and \eqref{eq:P-fac}, we obtain
$$g^N - q(g-1)^N = (g^N + \overline e_{N-1}g^{N-1} + \dots + \overline e_0)\cdot (\overline f_\ell g^\ell + \overline f_{\ell-1}g^{\ell-1} + \dots + \overline f_0)$$
(where the bar indicates the evaluation at 0). These substitutions make sense since we have established polynomiality of the coefficients in the previous paragraph. Consequently $$\overline f_0 = 1 - q$$ and $\overline f_m = 0$ for $m > 0$. Then $\overline e_i$ is the coefficient of $g^i$ in $$\overline f_0^{-1} (g^N - q(g-1)^N)=\frac{g^N - q(g-1)^N}{1-q},$$ which is clearly a rational function in $q$. It follows that $\overline {\mathsf R}$ is a rational function in $q$.

Next, we claim that the derivatives $$\frac{\partial \mathsf R}{\partial x_j}$$ are also given by rational functions in  $\{g_i\}$, $\{w_i\}$, $\{x_m\}$ and $q$, symmetric in the $\{g_i\}$'s. (We are viewing $g_i$ as functions of the independent variables $\{w_i\}$, $\{x_m\}$ and $q$.) Indeed, we have 
$$\frac{\partial \mathsf R}{\partial x_j} = \sum_{i=1}^{N} \frac{\partial \mathsf R}{\partial g_i} \cdot \frac{\partial g_i}{\partial x_j} + \partial_{j} \mathsf R.$$
(Here $\partial_{j}\mathsf R$ means to take the $x_j$-derivative treating the $\{g_i\}$ as constants.)
The second term $\partial_{j} \mathsf R$ is manifestly symmetric in the $g$'s since $\mathsf R$ is. Again because $\mathsf R$ is symmetric, we see that transposing $g_i$ and $g_k$ turns $$\frac{\partial \mathsf R}{\partial g_i} \text{ into  }\frac{\partial \mathsf R}{\partial g_k}.$$ For each fixed $i$, the derivative $\frac{\partial g_i}{\partial x_j}$ can be expressed as a rational function in the $g_i$ (but no other $g$'s), the $\{w_i\}$, $\{x_m\}$, and $q$ by implicit differentiation applied to $$ \mathsf X(g_i)=(-1)^{N} q.$$ Replacing $g_i$ with $g_k$ in this formula therefore yields the formula for $\frac{\partial g_k}{\partial x_j}$. The claim now follows since we sum over all $i$.

Inductively, it follows that all higher derivatives of $\mathsf R$ are rational functions of  $\{g_i\}$, $\{w_i\}$, $\{x_m\}$ and $q$, symmetric in the $\{g_i\}$. This completes the proof of Theorem \ref{t2}.
\qed

\subsection{Example}\label{examp} The proof of the Theorem gives an effective algorithm of computing the descendent series for $\beta=0$. We illustrate the case $$N=2,\,\, \alpha_1=\mathcal O_X, \,\, k_1=1,\,\, \ell=1.$$ Thus $$Z=\sum_{n=0}^{\infty} q^n\cdot  \int_{\left[\mathsf {Quot}_{X}(\mathbb C^2, n)\right]^{\text{vir}}} \mathsf{ch}_{1}(\mathcal O^{[n]}) \cdot \, c(T^{\text{vir}}\mathsf {Quot}).$$ When $N=1$, for the Hilbert scheme of points, the boundary insertion $c_1(\mathcal O^{[n]})$ plays an important role in the formalism of \cite {Le}. 

Setting $$W=\sum_{n=0}^{\infty} q^n\cdot  \int_{\left[\mathsf {Quot}_{X}(\mathbb C^2, n)\right]^{\text{vir}}} c_x(\mathcal O^{[n]}) \cdot \, c(T^{\text{vir}}\mathsf {Quot}),$$ we have $$W=\mathsf A^{K_X^2}$$ for some universal series $\mathsf A$. No other universal functions are needed in this case. Thus 
 $$Z=\frac{\partial W}{\partial x}\bigg\rvert_{x=0}=K_X^2\cdot \mathsf A^{K_X^2-1}\cdot \frac{\partial \mathsf A}{\partial x}\bigg\rvert_{x=0}.$$ We already calculated 
$$\mathsf A\bigg\rvert_{x=0} =\frac{(1-q)^2\cdot (1-6q+q^2)}{(1-4q)^2}$$ in Theorem 18 of \cite {OP1}. We furthermore claim \begin{equation}\label{de}\frac{\partial \mathsf A}{\partial x}\bigg\rvert_{x=0}=  \frac{2q^2 \cdot (1 - 12 q - 33 q^2 + 8 q^3)}{(1 - 4q)^3}.\end{equation} This follows by the proof of Theorem \ref{t2}. Indeed, we have $$\mathsf A^{-1}=\sum_{q=0}^{\infty} (-q)^n \int_{\mathsf {Quot}_{{\mathbb P}^1}(\mathbb C^2,n)} c_{x}(\OO^{[n]})  \cdot c(T\mathsf{Quot}_{\mathbb P^1})  \cdot c(\mathcal T_n)=-\frac{\Psi}{K}(g_1+w_1, g_2+w_2)$$ where $g_1, g_2$ solve the equation $$(g+w_1)\cdot(g+w_2)=q\cdot (1-g-w_1)\cdot (1-g-w_2)\cdot (1-xg).$$ The expressions for $\Psi, \Phi_i, K$ are explicitly given in the proof of the Theorem. Substituting and carrying out the implicit differentiation with respect to $x$, we arrive at expression \eqref{de} claimed above. 

\section{Descendent series for the Hilbert scheme: Theorem \ref{t3}}
\subsection{Descendents}
The argument of Theorem \ref{t1} extends to prove the more general descendent
claim of Theorem \ref{t3}. For each $K$-theory class $\alpha$ on $X$,
we have defined $$\alpha^{[n]}={\mathbf R}\pi_{\star}\left(\mathcal Q\otimes p^{\star} \alpha\right)$$ on $\mathsf {Quot}_{X}(\mathbb C^N, \beta, n)\, .$
The descendent series is given by  \begin{multline*}Z_{X, \, N,\,\beta}(\alpha_1, \ldots, \alpha_{\ell}\, |\, k_1, \ldots, k_{\ell})=\\
  \sum_{n\in \mathbb{Z}}  q^n\cdot \int_{\left[\mathsf {Quot}_X(\mathbb C^N, \beta, n)\right]^{\text{vir}}} \mathsf {ch}_{k_1} (\alpha_1^{[n]}) \cdots \mathsf{ch}_{k_{\ell}} (\alpha_{\ell}^{[n]})\cdot c(T^{\text{vir}} \mathsf {Quot})\, .\end{multline*} 

\noindent To establish Theorem \ref{t3}, we set $N=1$, and show that the series $$\generall\in\mathbb Q((q))$$ is the Laurent expansion of a rational function for a nonsingular projective simply connected surface $X$. 

\subsection {Proof of Theorem \ref{t3}} 

\subsubsection{Hilbert schemes of points}
We use again the isomorphism $$\mathsf {Quot}_X(\mathbb C^1, \beta, n)\simeq X^{[m]}\times \mathbb P\, ,\ \ \
m=n+\frac{\beta(\beta+K_X)}{2}$$ where $\mathbb P$ denotes the linear system $|\beta|$. We will study the series
\begin{equation}\label{zt2}Z=\sum_{n\in \mathbb{Z}} q^n\cdot \int_{X^{[m]}\times \mathbb P}\mathsf {ch}_{k_1} (\alpha_1^{[n]}) \cdots \mathsf{ch}_{k_{\ell}} (\alpha_{\ell}^{[n]})\cdot \mathsf {e}(\mathsf{Obs})\cdot \frac{c(TX^{[m]}) c(T\mathbb P)}{c(\mathsf{Obs})}\, .\end{equation}

We identify the tautological structures appearing in \eqref{zt2}. The universal quotient over $\mathsf {Quot}_
X(\mathbb C^1, \beta, n) \times X$ can be expressed in $K$-theory as $$\mathcal Q=\mathcal O - \mathcal I_{\mathcal W} \otimes \mathcal O(-\beta)\otimes \mathcal L^{-1}=\mathcal O  - \mathcal O(-\beta)\otimes\mathcal L^{-1} + \mathcal O_{\mathcal W}\otimes \mathcal O(-\beta)\otimes\mathcal L^{-1}$$ where
$\mathcal W$ denotes the universal subscheme of $X^{[m]}\times X,$
and $$\mathcal L=\mathcal O_{\mathbb P}(1)\to \mathbb P$$
denotes
the tautological bundle. As a result $$\alpha^{[n]}_{\mathsf {Quot}}=H^{\bullet}(\alpha) \otimes \mathcal O-H^{\bullet}(\widetilde \alpha)\otimes \mathcal L^{-1}+ \widetilde \alpha^{[m]}_{\mathsf {Hilb}}\otimes \mathcal L^{-1}$$
where $\widetilde \alpha=\alpha \otimes \mathcal O(-\beta)$.
We have
indicated by subscripts the locations of the tautological constructions.
Let $\zeta=c_1(\mathcal L).$ Thus $$\mathsf {ch} (\alpha^{[n]}_{\mathsf {Quot}})=\chi(\alpha)-\chi(\widetilde \alpha)\cdot e^{-\zeta}+ \mathsf {ch}(\widetilde \alpha^{[m]}_{\mathsf {Hilb}})\cdot e^{-\zeta}$$ which, in fixed degree $k>0$, becomes $$\mathsf {ch}_{k} (\alpha^{[n]}_{\mathsf {Quot}})=-\chi(\widetilde \alpha)\cdot\frac{(-\zeta)^{k}}{k!} + \sum_{j=0}^{k} \mathsf {ch}_{j}(\widetilde \alpha^{[m]}_{\mathsf{Hilb}}) \cdot \frac{(-\zeta)^{k-j}}{(k-j)!}\, .$$
After multiplying out the different Chern characters appearing in \eqref{zt2},
we are led to expressions of the form \begin{equation}\label{newx}
\sum_{m=0}^{\infty} q^m \cdot \int_{X^{[m]}\times \mathbb P}\zeta^{k}\cdot \mathsf {ch}_{k_1} (\alpha_1^{[m]}) \cdots \mathsf{ch}_{k_{\ell}} (\alpha_{\ell}^{[m]})\cdot \mathsf {e}(\mathsf{Obs})\cdot \frac{c(TX^{[m]}) c(T\mathbb P)}{c(\mathsf{Obs})}\, .
\end{equation} Here, we have changed notation by
removing the tilde's from the $\alpha$'s and relabeling indices.{\footnote{The
    overall $q$ shift does not affect rationality.}}
Unless specified otherwise, all tautological structures $\alpha^{[m]}$ are from now on understood to be over the Hilbert scheme of points $X^{[m]}$.

We will consider two cases depending upon
the geometric genus of the simply connected surface $X$. Furthermore, when the genus is positive, we first discuss surfaces which are minimal, and then consider their blowups. 

\subsubsection{Minimal surfaces with $p_g>0$} \label{minpg} Assume that $X$ is simply connected minimal surface. Then $X$ is either a $K3$ surface, an elliptic surface, or a surface of general type. \vskip.1in

\noindent $\bullet$ For $K3$ surfaces, the virtual fundamental class vanishes due to the presence of a trivial factor in the obstruction bundle, unless $\beta=m=0$ \cite{MOP}. There is nothing to prove in the $K3$ case. \vskip.1in

\noindent $\bullet$ If $X$ minimal of general type, the virtual fundamental class of $\mathsf {Quot}_{X}(\mathbb C^1, n, \beta)$ was shown to vanish in \cite [Section $5.3.3$]{OP1}, unless 
\begin{itemize}
\item [(i)] $\beta=0$ or 
\item [(ii)] $\beta=K_X$ and $m=0$.
\end{itemize} 
\noindent There is nothing to prove in case (ii). When $\beta=0$, we can use Theorem \ref{t2} or alternatively, we can argue as follows. We have $$\mathsf {Obs}=\left((K_X)^{[m]}\right)^{\vee},$$ see for instance \eqref{obun}. The series \eqref{newx} becomes $$Z=\sum_{m=0}^{\infty} q^m \int_{X^{[m]}}\mathsf {ch}_{k_1} (\alpha_1^{[m]}) \cdots \mathsf{ch}_{k_{\ell}} (\alpha_{\ell}^{[m]})\cdot \mathsf e\left(\left(K_X^{[m]}\right)^{\vee}\right) \cdot \frac{c(TX^{[m]})}{c \left(\left(K_X^{[m]}\right)^{\vee}\right)}\, .$$ We conclude by Proposition \ref{rat} below. 
\vskip.1in

\noindent $\bullet$ Let  $X\to C$ be a minimal elliptic surface with $p_g>0$. Since $X$ is simply connected, we must have $C=\mathbb P^1$, by \cite[Lemma VII.14]{Fr}. We first argue that $\beta$ must be a multiple of the fiber. Note that $\beta$ must be effective for the Quot scheme to be nonempty. Furthermore, the expression for the obstruction bundle \eqref{obun}, $$\mathsf{Obs} = (H^1(M)-H^0(M))\otimes \mathcal L+\left(M^{[m]}\right)^{\vee}\otimes \mathcal L+ \mathbb C^{p_g}\, ,$$ shows that the virtual fundamental class vanishes if $$H^0(M)=0\iff H^0(K_X-\beta)= 0$$ due to the presence of the trivial factor. We may therefore
assume $K_X-\beta$ is effective. Since $X$ is minimal, we find $$K_X=(p_g-1)f\, .$$ Since $$(p_g-1)f=\beta+(K_X-\beta)$$ is an effective decomposition,  $\beta$ must be supported on fibers. By Zariski's Lemma, $\beta^2\leq 0$ and $\beta\cdot K_X=0$. If $\beta^2<0$ then \begin{equation}\label{inq}\beta\cdot (\beta-K_X)<0\, .\end{equation} When inequality \eqref{inq}
is satisfied, the proof of Proposition $22$ of \cite {OP1} shows that the virtual fundamental class vanishes. Proposition 22 of \cite{OP1} is stated for surfaces of general type, but the same argument applies here as
well.\footnote{This can also be seen via \eqref{wl} since $\text{Hilb}_{\beta}$ has negative virtual dimension.} Thus $$\beta^2=0\, ,$$ so by Zariski's Lemma $\beta=af$ for $0\leq a\leq p_g-1.$

Recording that $$\chi(\mathcal O(af))=1+p_g\, ,\ \
h^0(\mathcal O(af))=a+1\, ,$$ we find that \eqref{obun} becomes $$\mathsf {Obs}=-\mathbb C^{p_g-a}\otimes \mathcal L +\left(M^{[m]}\right)^{\vee} \otimes \mathcal L+\mathbb C^{p_g}$$ over $X^{[m]}\times \mathbb P^{a}\,.$
We then obtain \begin{eqnarray*}\mathsf e(\mathsf {Obs})&=&\left[c(\mathcal L)^{a-p_g}\cdot c\left(\left(M^{[m]}\right)^{\vee} \otimes \mathcal L\right)\right]_{(a+m)}\\
&=&\left[(1+\zeta)^{a-p_g} \cdot \sum_{k=0}^{m} (1+\zeta)^{k} c_{m-k} \left(\left(M^{[m]}\right)^{\vee}\right)\right]_{(a+m)}.\end{eqnarray*} The exponents of the hyperplane class $\zeta$ over $\mathbb P^a$ must be bounded by $a$. Thus, for degree reasons, the only contribution reaching the necessary degree $a+m$ occurs for $k=0$ and in this case \begin{equation}\label{inqq}\mathsf e(\mathsf {Obs})=\left[(1+\zeta)^{a-p_g}\right]_{(a)}\times c_m\left(\left(M^{[m]}\right)^{\vee}\right)=\binom{a-p_g}{a} \left[\text{pt}\right] \times \mathsf e\left(\left(M^{[m]}\right)^{\vee}\right)\end{equation} over $\mathbb P^a\times X^{[m]}.$ 

As a result of the above calculation, the series \eqref{newx}, $$\sum_{m=0}^{\infty} q^m \int_{X^{[m]}\times \mathbb P^a}\zeta^{k}\cdot \mathsf {ch}_{k_1} (\alpha_1^{[m]}) \cdots \mathsf{ch}_{k_{\ell}} (\alpha_{\ell}^{[m]})\cdot \mathsf {e}(\mathsf{Obs})\cdot \frac{c(TX^{[m]}) c(T\mathbb P)}{c(\mathsf{Obs})}\, ,$$ vanishes for $k>0$. For $k=0$, the expression simplifies to $$\binom{a-p_g}{a} \cdot \sum_{m=0}^{\infty} q^m \int_{X^{[m]}} \mathsf {ch}_{k_1} (\alpha_1^{[m]}) \cdots \mathsf{ch}_{k_{\ell}} (\alpha_{\ell}^{[m]})\cdot \mathsf e\left(\left(M^{[m]}\right)^{\vee}\right) \cdot \frac{c(TX^{[m]})}{c\left(\left(M^{[m]}\right)^{\vee}\right)}\, .$$ Proposition \ref{rat} below completes the argument. 

Theorem \ref{t3} is established for all simply connected minimal surfaces with $p_g>0$. 

\subsubsection{Further descendent rationality} \label{r1} We prove here
the following result that was used in Subsection \ref{minpg}. 
\begin{proposition}\label{rat} The generating series $$\sum_{n=0}^{\infty} q^n\cdot
  \int_{X^{[n]}}\mathsf {ch}_{k_1} (\alpha_1^{[n]}) \cdots \mathsf{ch}_{k_{\ell}} (\alpha_{\ell}^{[n]})\cdot \mathsf e\left(\left(M^{[n]}\right)^{\vee}\right) \cdot \frac{c(TX^{[n]})}{c \left(\left(M^{[n]}\right)^{\vee}\right)}$$ is
  a rational function in $q$ for all pairs $(X, M).$ 
\end{proposition} 

\proof The proof is similar to that of Theorem \ref{t2}, using the methods developed in \cite {OP1}. Just as in Theorem \ref{t2}, we may assume that $$\text{rank }\alpha_i=1 \text{ for all }1\leq i\leq \ell.$$ Expressing the Chern character in terms of Chern classes, it suffices to show that the series
$$Z=\sum_{n=0}^{\infty} q^n \int_{X^{[n]}} c_{k_1}(\alpha_1^{[n]}) \dots c_{k_\ell}(\alpha_\ell^{[n]})\cdot \mathsf e\left(\left(M^{[n]}\right)^{\vee}\right) \cdot \frac{c(TX^{[n]})}{c \left(\left(M^{[n]}\right)^{\vee}\right)}$$ is a rational function in $q$. 

Let $x_1, \ldots, x_{\ell}$ be formal variables. Write $$c_x=1+x c_1+ x^2 c_2 + \ldots$$ for the total Chern class, and set $$W=\sum_{n=0}^{\infty} q^n \int_{X^{[n]}} c_{x_1}(\alpha_1^{[n]}) \dots c_{x_\ell}(\alpha_\ell^{[n]})\cdot \mathsf e\left(\left(M^{[n]}\right)^{\vee}\right) \cdot \frac{c(TX^{[n]})}{c \left(\left(M^{[n]}\right)^{\vee}\right)}.$$ The series $Z$ is found by extracting the coefficient of $x_1^{k_1}\cdots x_{\ell}^{k_{\ell}}$ in $W$: $$Z=\frac{1}{k_1!}\cdots \frac{1}{k_\ell!}\cdot \frac{\partial^{k_1}}{\partial^{k_1}x_1} \ldots \frac{\partial^{k_\ell}}{\partial^{k_\ell}x_\ell} W\bigg\rvert_{x_1=\ldots=x_{\ell}=0}.$$ Now, invoking the universality and multiplicativity results of \cite {EGL}, we find the factorization 
$$W = {\mathsf A}^{K_X^2} \cdot {\mathsf B}^{\chi(X)} \cdot {\mathsf C}^{M^2} \cdot {\mathsf D}^{M.K_X}\cdot \prod_{i=1}^{\ell} {\mathsf E}_i^{c_1(\alpha_i).K_X} \cdot \mathsf F_i^{c_2(\alpha_i)}\cdot \mathsf G_i^{c_1(\alpha_i)\cdot M} \cdot \prod_{1\leq i\leq j\leq \ell} \mathsf H_{ij}^{c_1(\alpha_i)\cdot c_1(\alpha_j)}$$
in terms of universal series that depend on $q$ and $x_i$. To find these series, we can pick convenient geometries. We may assume $M$ is sufficiently positive, so that there exists $C$ a nonsingular connected curve in the linear system $|M|$. As explained in \cite {OP1}, we have $$\mathsf e\left(\left(M^{[n]}\right)^{\vee}\right)\cap \left[X^{[n]}\right] =(-1)^{n} j_{\star}\left[C^{[n]}\right]$$ where $$j:C^{[n]}\hookrightarrow X^{[n]}$$ is the natural inclusion. By equation (33) of \cite {OP1}, we furthermore have $$j^{\star} c(TX^{[n]})=c\left(\left(K_C^{[n]}\right)^{\vee}\right) \cdot c(M^{[n]}).$$ Thus $$W=\sum_{n=0}^{\infty} (-q)^n \int_{C^{[n]}} c_{x_1} (\alpha_1^{[n]}) \cdots c_{x_{\ell}} (\alpha_{\ell}^{[n]}) \cdot \frac{c\left(\left(K_C^{[n]}\right)^{\vee}\right) \cdot c(M^{[n]})}{c\left(\left(M^{[n]}\right)^{\vee}\right)}$$ or equivalently, in terms of Segre classes
\begin{multline*}\sum_{n=0}^{\infty} (-q)^n \cdot \int_{C^{[n]}} s_{x_1} \left((-\alpha_1)^{[n]} \right)\cdots s_{x_{\ell}} \left((-\alpha_\ell)^{[n]}\right) \\ \cdot s_1\left((-M)^{[n]}\right)\cdot s_{-1}\left((-K_C)^{[n]}\right) \cdot s_{-1}(M^{[n]})\, .
\end{multline*} Here, $\alpha_1, \cdots, \alpha_\ell, M$ are understood to be restricted from the surface $X$ to the curve $C$.

Using Theorem 3 of \cite {OP1}, the last expression can be evaluated in closed form. Under the change of variables $$-q=\frac{t}{(1-t)(1-x_1t)\cdots (1-x_{\ell}t)}\, ,$$ we have \begin{equation}\label{ew}W(q, x_1, \ldots, x_\ell)
  =\prod_{i=1}^{\ell} \left(1-x_it\right)^{-c_1(\alpha_i)\cdot M} \cdot (1-t)^{-M^2} \cdot (1+t)^{-M\cdot K_X} \cdot \mathsf a^{K_X^2}\end{equation} with{\footnote{We can show $$\mathsf a=1-e_2t^2+2e_3t^3-3e_4t^4+\ldots $$ where $e_i$ are the elementary symmetric functions in $1, x_1, \ldots, x_{\ell}$.
    We do not explain the latter formula for $\mathsf a$ since it will
    not be used here.}}
$$\mathsf a=-\left(\frac{q}{t}\right)^{-2} \cdot \frac{dq}{dt}\, .$$
It follows that $$\mathsf {B}=1,\quad \mathsf{E_i}=1,\quad \mathsf{F_i}=1, \quad \mathsf{H}_{ij}=1.$$ Furthermore, by universality, expression \eqref{ew} for $W$ holds for all geometries $(X, M, \alpha_1, \ldots, \alpha_{\ell})$, not only for those for which $M$ is sufficiently positive. 

Identity \eqref{ew} for the values $x_1=\ldots=x_{\ell}=0$ implies
$$W(q, 0, \ldots, 0)=(1-t)^{-M^2} \cdot (1+t)^{-M\cdot K_X}$$ with $$q=-\frac{t}{1-t}\iff t=-\frac{q}{1-q}\, .$$ Evidently $W(q, 0, \ldots, 0)$ is a rational function of $t$ and hence also of $q$. In fact, the expression we have obtained, \begin{equation}\label{wq0} W(q, 0, \ldots, 0)=(1-q)^{M^2}\cdot \left(\frac{1-q}{1-2q}\right)^{M\cdot K_X}\, ,\end{equation} is Corollary 38 of \cite {OP1}. 

However, we can now also go further. We address all derivatives of $W$ with respect to $x_i$ via \eqref{ew}, as needed to complete the proof of Proposition \ref{rat}. Clearly, the derivatives $$\frac{\partial^{k_1}}{\partial^{k_1}x_1} \ldots \frac{\partial^{k_\ell}}{\partial^{k_\ell}x_\ell}\bigg\rvert_{x_1=\ldots=x_{\ell}=0} \,\text{RHS of }\eqref{ew}$$ are rational functions in $t$. For the left hand side of \eqref{ew}, we apply the chain rule repeatedly using that $$\frac{\partial q}{\partial x_i}\bigg\rvert_{x_1=\ldots=x_{\ell}=0}=-\frac{t^2}{1-t},\,\,,\,\,\, q\bigg\rvert_{x_1=\ldots=x_{\ell}=0}=-\frac{t}{1-t}.$$ For instance, the $x_1$ derivative equals $$\frac{d}{d x_1} W(q, x_1, \ldots, x_{\ell})\bigg\rvert_{x_1=\ldots=x_\ell=0}=\frac{\partial W}{\partial q}\bigg\rvert_ {\stackrel{q=-\frac{t}{1-t}}{{x_1=\ldots=x_{\ell}=0}}}\cdot \frac{-t^2}{1-t} + \frac{\partial W}{\partial x_1}\bigg\rvert_ {\stackrel{q=-\frac{t}{1-t}}{x_1=\ldots=x_\ell=0}}.$$ We argued the left hand side is rational in $t$. Since $\frac{\partial W}{\partial q} \left(q, 0, \ldots, 0\right)$ is rational in $t$ by \eqref{wq0}, we conclude that the same is true about the last term $$\frac{\partial W}{\partial x_1}\bigg\rvert_ {\stackrel{q=-\frac{t}{1-t}}{x_1=\ldots=x_\ell=0}}.$$ Equivalently, $$\frac{\partial W}{\partial x_1}\bigg\rvert_{x_1=\ldots=x_{\ell}=0}$$ is rational in $q$. 
Rationality of the higher order derivatives follows inductively. The proof of
Proposition \ref{rat} is complete. 
\qed
\vskip.1in

\subsubsection{Example}\label{r8} We illustrate Proposition \ref{rat} with the computation of the series\footnote{The same calculation can also be carried out using Theorem \ref{t2}. } $$Z=\sum_{n=0}^{\infty} q^n\cdot  \int_{\left[X^{[n]}\right]^{\text{vir}}} \mathsf{ch}_{1}(\mathcal O^{[n]}) \cdot \, c(T^{\text{vir}}X^{[n]})\, .$$ 
This expression is of the form considered in the Proposition. Indeed, by equation $(33)$ of \cite {OP1} we have $$\left[X^{[n]}\right]^{\text{vir}}=\mathsf e\left(\left(M^{[n]}\right)^{\vee}\right),\quad T^{\text{vir}}X^{[n]}=TX^{[n]}-\left(M^{[n]}\right)^{\vee}$$ for $M=K_X$.
Write 
$$W=\sum_{n=0}^{\infty} q^n\cdot   \int_{\left[X^{[n]}\right]^{\text{vir}}}  c_x(\mathcal O^{[n]}) \cdot \,c(T^{\text{vir}}X^{[n]})\implies Z=\frac{\partial W}{\partial x}\bigg\rvert_{x=0}.$$ Equation \eqref{ew} for $M=K_X$ expresses the answer in terms of a single universal series $$W=\mathsf A^{K_X^2}\,$$ where $$\mathsf A(q, x)=\frac{1-t^2x}{1-t^2}, \,\, \text { for } q=-\frac{t}{(1-t)(1-tx)}\, .$$ Thus $$Z=K_X^2\cdot \mathsf A^{K_X^2-1}\cdot \frac{\partial \mathsf A}{\partial x}\bigg\rvert_{x=0}\, .$$ By direct computation
$$\mathsf A\bigg\rvert_{x=0} =\frac{(1-q)^2}{1-2q},\,\,\,\,\,\,\, \frac{\partial \mathsf A}{\partial x}\bigg\rvert_{x=0}= - \frac{q^2\cdot (1 - 4 q)}{(1-2q)^2}\, .$$ In fact, it can easily be shown that all derivatives take the form $$\frac{\partial^k \mathsf A}{\partial x^k}\bigg\rvert_{x=0}=\frac{\mathsf P_k(q)}{(1-q)^{2k-2}\cdot (1-2q)^{k+1}}$$ for some polynomials $\mathsf P_k$. 
The denominators $1-q$ and $1-2q$ are consistent with the proof of Theorem \ref{t1}. 
\vskip.1in
\subsubsection{Non-minimal surfaces with $p_g>0$}
\label{nms}
For non-minimal surfaces, we prove Theorem \ref{t3} using the calculations of
Section \ref{minpg}, combined with an observation that we learned from Woonam Lim. Specifically, in the next paragraph, we will explain a special case of Lemma 2 of \cite {L} in our simpler setting. 
The argument rests on a deeper connection with Seiberg-Witten theory and the notion of simple type. 

Let $\widetilde X$ denote the blowup of a nonsingular projective
simply connected
surface $X$  with exceptional divisor $E$. For each curve class $\beta$ on $X$ and each integer $k$, consider the class $$\widetilde \beta=\beta+k E$$ on $\widetilde X$. We assume
\begin{equation} \label{pp22}
  \left[\mathsf {Quot}_{X}(\mathbb C^1, \beta, n)\right]^{\text{vir}}\neq 0 \text{ for some } n\implies \beta\cdot (\beta-K_X)=0\, .
\end{equation}
Assumption \eqref{pp22} is satisfied for all three classes of minimal surfaces
 considered in Section \ref{minpg} as the reader can immediately verify. We claim that the same holds true on $\widetilde X$: 
$$\left[\mathsf {Quot}_{\widetilde X}(\mathbb C^1, \widetilde \beta, n)\right]^{\text{vir}}\neq 0 \text { for some }n\implies \widetilde \beta\cdot (\widetilde \beta-K_{\widetilde X})=0\, .$$ 
By direct calculation, $$\widetilde \beta\cdot (\widetilde \beta-K_{\widetilde X})\ = \ \beta\cdot (\beta-K_X)-k(k-1)\ =\ -k(k-1)\leq 0\, .$$
If the inequality is strict $$\widetilde \beta\cdot (\widetilde \beta-K_{\widetilde X})<0\, ,$$ the virtual fundamental class vanishes by the proof of Proposition 22 of \cite {OP1} (as already used in equation \eqref{inq}).
We must therefore
have $\widetilde \beta \cdot (\widetilde \beta-K_{\widetilde X})=0$.  

Applying the argument inductively to a sequence of blowups, we see that if $X$ is a possibly non-minimal surface with $p_g>0$, non-zero invariants only arise if \begin{equation}\label{swc}\beta\cdot (\beta-K_X)=0\, .\end{equation}

The latter condition can be used to explicitly calculate the virtual fundamental class. Indeed, thanks to \eqref{swc}, and recalling the obstruction bundle from equation \eqref{obun}, we have $$\text{rank }\mathsf {Obs}=m+h^0(\beta)-1\,
.$$ We now use the same reasoning that led to \eqref{inqq}. For the current numerics, we similarly compute over $X^{[m]}\times \mathbb P$:
\begin{eqnarray*}\mathsf e({\mathsf {Obs}})&=&\left[c(\mathcal L)^{h^1(\beta)-h^2(\beta)}\cdot c\left(\left(M^{[m]}\right)^{\vee} \otimes \mathcal L\right)\right]_{m+h^0(\beta)-1}
  \\&=&\left[(1+\zeta)^{h^1(\beta)-h^2(\beta)}\cdot \sum_{k=0}^{m} (1+\zeta)^{k} c_{m-k} \left(\left(M^{[m]}\right)^{\vee}\right)\right]_{m+h^0(\beta)-1}\\&=&\binom{h^1(\beta)-h^2(\beta)}{h^0(\beta)-1} \cdot \mathsf e\left( \left(M^{[m]}\right)^{\vee}\right)\times \left[\text{pt}\right].\end{eqnarray*} The argument then
is completed in the same fashion as for elliptic surfaces in Section \ref{minpg}  by invoking Proposition \ref{rat}. 

\subsubsection {Surfaces with $p_g=0$.} We establish Theorem \ref{t3} for surfaces with $p_g=0$. We follow here
the proof in Section \ref{rational} closely. 
We have $$\mathsf {Obs}=H^1(M)^{\vee}\otimes \mathcal L + \left(M^{[n]}\right)^{\vee}\otimes \mathcal L\, .$$ By \eqref{newx}, we examine expressions of the form $$\sum_{n=0}^{\infty} q^n \int_{X^{[n]}\times \mathbb P}\zeta^{k+h^1(\beta)}\cdot \mathsf {ch}_{k_1} (\alpha_1^{[n]}) \cdots \mathsf{ch}_{k_{\ell}} (\alpha_{\ell}^{[n]})\cdot \mathsf e\left(\mathcal L\otimes \left(M^{[n]}\right)^{\vee}\right) \cdot \frac{c(TX^{[n]}) \cdot c(\mathcal L)^{\chi(\beta)}}{c\left(\mathcal L\otimes \left(M^{[n]}\right)^{\vee}\right)}\, .$$ Expanding the terms that involve $\mathcal L$ into powers of $\zeta=c_1(\mathcal L)$ as in Proposition \ref{p1a}, we obtain $$\sum_{n=0}^{\infty}q^n\int_{X^{[n]}\times \mathbb P}\zeta^{k+h^1(\beta)}\cdot \mathsf {ch}_{k_1} (\alpha_1^{[n]}) \cdots \mathsf{ch}_{k_{\ell}} (\alpha_{\ell}^{[n]})\cdot \frac{c(TX^{[n]})}{c\left( \left(M^{[n]}\right)^{\vee}\right)}$$ $$\cdot \left(\sum_{a=0}^{n} \zeta^a \cdot c_{n-a} \left(\left(M^{[n]}\right)^{\vee}\right)\right)\cdot \left(\sum_{j=0}^{\infty} (-1)^{j} \zeta^{j} H_j \right)\cdot (1+\zeta)^{\chi(\beta)} .$$ Integrating out the powers of $\zeta$, we equivalently prove the rationality of 
\begin{equation}\label{aa}\sum_{n=0}^{\infty} q^n\int_{X^{[n]}} \mathsf {ch}_{k_1} (\alpha_1^{[n]}) \cdots \mathsf{ch}_{k_{\ell}} (\alpha_{\ell}^{[n]})\cdot c_{n-a}\left(\left(M^{[n]}\right)^{\vee}\right) \cdot c(TX^{[n]}) \cdot \frac{P\left(\left(M^{[n]}\right)^{\vee},B\right)}{c\left(\left(M^{[n]}\right)^{\vee}\right)}\, ,\end{equation} for fixed tuples $(a, B, k_1, \ldots, k_{\ell}, \alpha_1, \ldots, \alpha_{\ell})$. Following the proof of Proposition \ref{p1a}, we will establish first universality and then rationality for sufficiently positive geometries. 

For universality, we first turn all Chern characters into universal expressions in the Chern classes: 
$$\sum_{n=0}^{\infty} q^n\int_{X^{[n]}} c_{k_1} (\alpha_1^{[n]}) \cdots c_{k_{\ell}} (\alpha_{\ell}^{[n]})\cdot c_{n-a}\left(\left(M^{[n]}\right)^{\vee}\right) \cdot c(TX^{[n]}) \cdot \frac{P\left(\left(M^{[n]}\right)^{\vee},B\right)}{c\left(\left(M^{[n]}\right)^{\vee}\right)}\, .$$
We introduce formal variables $x_1, \ldots, x_{\ell}$, 
and form the generating series 
\begin{multline*}
  Y_{X, M}^{(p)}=\sum_{B=(b_1, \ldots, b_p)} \frac{z_1^{b_1}}{b_1!} \cdots \frac{z_{p}^{b_{p}}}{b_{p}!} \sum_{n\geq 0} \sum_{a\geq 0}q^n t^a \cdot \int_{X^{[n]}} c_{x_1} (\alpha_1^{[n]}) \cdots c_{x_{\ell}} (\alpha_{\ell}^{[n]})\ c_{n-a}\left(\left(M^{[n]}\right)^{\vee}\right)\\ \cdot\,  c(TX^{[n]}) \cdot \frac{P\left(\left(M^{[n]}\right)^{\vee},B\right)}{c\left(\left(M^{[n]}\right)^{\vee}\right)}
  \, .
\end{multline*}
The length of $B$ equals the superscript $p$ appearing on the left hand side.
We must
extract $$\text{Coefficient of}\ \ x_1^{k_1}\cdots x_{\ell}^{k_{\ell}}\cdot  \frac{z_1^{b_1}}{b_1!} \ \cdots\ \frac{z_{p}^{b_{p}}}{b_{p}!} \cdot t^a
\ \ \text {in}\ \ Y^{(p)}_{X, M}\, .$$
As in Section \ref{uni}, $Y_{X, M}^{(p)}$ is multiplicative and can be factored in terms of several universal power series. It suffices therefore to
establish rationality (of the correct coefficient) for special geometries. 

Returning to expression \eqref{aa}, we pick a sufficiently positive $M$, and represent $c_{n-a}\left(M^{[n]}\right)$ by the relative Hilbert scheme $$(\mathcal C/\mathcal B)^{[n]}\to \mathcal B$$ of a linear system $|V|\subset |M|$ as in Section \ref{uni2}. By the arguments of the same Section, it suffices to consider expressions of the form $$\sum_{n=0}^{\infty} q^n\int_{(\mathcal C/\mathcal B)^{[n]}} \mathsf{ch}_{k_1} (j^{\star}\alpha_1^{[n]}) \cdots \mathsf {ch}_{k_{\ell}} (j^{\star}\alpha_{\ell}^{[n]}) \cdot c(\gamma+\alpha_n\cdot \zeta_n^{-1}+\beta_n\cdot \zeta_n)\cdot P\left(j^{\star} \left(M^{[n]}\right)^{\vee}, B\right),$$ where,
as before,
$$j:(\mathcal C/\mathcal B)^{[n]}\to X^{[n]}\, .$$ Let $\mu$ denote one of the classes $\alpha_1, \ldots, \alpha_{\ell}.$ Invoking \eqref{comput}, we have $$j^{\star} \mu^{[n]}=\overline {\mu}-\mu_n\cdot \zeta_n^{-1}$$ and
hence \begin{equation}\label{sss}\mathsf {ch}_{k }\, (j^{\star}\mu^{[n]})=\mathsf{ch}_k (\overline {\mu})-\sum_{i=0}^{k} \frac{(-1)^{k-i}}{(k-i)!}\cdot \mathsf{ch}_i(\mu_n)\cdot c_1(\zeta_n)^{k-i} .\end{equation}
Following the derivation of equation \eqref{eee}, we obtain  
\begin{equation}\label{eeee}\int_{\mathbb P_n} c_1(\zeta_n)^{s} \cdot \rho_n\cdot \lambda\cdot c(\alpha_n\cdot \zeta_n^{-1}+\left(\beta_n\right)^{\vee}\cdot \zeta_n)\cdot P_{b_1}\left((-M)_n^{\vee}\cdot \zeta_n\right)\cdots P_{b_m}\left((-M)_n^{\vee}\cdot \zeta_n\right).\end{equation} 
\noindent Compared to \eqref{eee}, the extra terms are $c_1(\zeta_n)^{s}$ and the class $\rho_n$ which is a universal polynomial in the Chern classes $$c_{i}(\mu_n)$$ where $\mu$ is one of the classes $\alpha_1, \ldots, \alpha_\ell$. These extra terms arise from the product expansion $$\mathsf{ch}_{k_1} (j^{\star}\alpha_1^{[n]}) \cdots \mathsf {ch}_{k_{\ell}} (j^{\star}\alpha_{\ell}^{[n]})$$ using \eqref{sss}. Crucially for us, $s$ and the $i$'s are bounded from
above by an expression that depends on $k_1, \ldots, k_\ell$. Thus they are independent of $n$. 

The rest of the argument is as in Sections \ref{uni2} and \ref{beee}: we expand all expressions in powers of $c_1(\zeta_n)$ and integrate over the fibers of
$$\mathbb P_n\to \mathfrak M\, .$$
Keeping track of the numerical modifications is not difficult. The powers $c_1(\zeta_n)^{s}$ affect the indices of various sums defining the prefactors $\sigma^{(n)}$, see for instance \eqref{vwvw}. Since $s$ is fixed independently of $n$, the conclusions of Lemma \ref{l3} still hold. Furthermore, Lemma \ref{l4} can be applied to each of the additional terms $c_j(\mu_n)$ for $\mu$ being one of $\alpha_1, \ldots, \alpha_\ell.$ In the end, \eqref{eeee} is still an expression of the form \eqref{special}. Rationality is therefore established. \qed

\end{document}